\documentclass[11pt]{article}
\usepackage[a4paper,margin=31mm]{geometry}
\usepackage{amsmath,amssymb,amsthm,mathtools,mathrsfs}
\usepackage{booktabs,array}
\usepackage{enumitem}
\usepackage{graphicx}
\usepackage{placeins}
\usepackage[hidelinks]{hyperref}
\hypersetup{
  pdftitle={Convergence analysis of physics-informed neural networks for subdiffusion},
  pdfauthor={Bangti Jin, Qiqi Rao, Zhi Zhou}
}
\usepackage{microtype}

\numberwithin{figure}{section}
\numberwithin{equation}{section}

\makeatletter\def\@captype{table}\makeatother
\usepackage{caption}
\date{}
\allowdisplaybreaks

\newtheorem{theorem}{\bf Theorem}[section]
\newtheorem{lemma}{\bf Lemma}[section]
\newtheorem{corollary}{\bf Corollary}[section]
\newtheorem{proposition}{\bf Proposition}[section]

\newtheorem{assumption}{\bf Assumption}[section]

\def\d{{\rm d}}

\def\Omega{\varOmega}

\title{Convergence analysis of physics-informed neural networks for subdiffusion\thanks{The work of B. J. is supported by Hong Kong RGC General Research Fund (Project
14306423 and Project 14306824), ANR / Hong Kong RGC Joint Research Scheme (A-CUHK402/24) and a start-up fund from The Chinese University of Hong Kong. The work of Z. Z. is supported by by National Natural
Science Foundation of China (Project 12422117 and Project 12426312) and Hong Kong Research
Grants Council (15302323), and an internal grant of Hong Kong Polytechnic
University (Project ID: P0053938, Work Programme: 4-ZZVA).}}

\author{Bangti Jin\thanks{Department of Mathematics, The Chinese University of Hong Kong, Shatin, New Territories, Hong Kong, P.R. China (\texttt{b.jin@cuhk.edu.hk, qqrao@math.cuhk.edu.hk})}
\and Qiqi Rao\footnotemark[2] 
\and Zhi Zhou\thanks{Department of Applied Mathematics, The Hong Kong Polytechnic University, Kowloon, Hong Kong, China. (\texttt{zhizhou@polyu.edu.hk})}
}

\begin{document}

\maketitle

\setlength\abovedisplayskip{4.5pt}
\setlength\belowdisplayskip{4.5pt}

\begin{abstract}
Physics-informed neural networks (PINNs) have emerged as mesh-free solvers for partial differential equations, yet their rigorous convergence theory for time-fractional problems remains limited. We present a convergence analysis of PINNs for linear subdiffusion equations and for the associated backward problem of recovering the initial data from terminal observations. To address the nonlocality and weakly singular kernels of fractional derivatives, we introduce a discretization-based PINN framework built on some popular time-stepping schemes satisfying discrete maximal regularity, including backward Euler, BDF2, L1, and fractional Crank--Nicolson schemes. This discretize-then-train approach avoids the exact evaluation of fractional derivatives at training points while preserving the stability structure needed for convergence analysis. For the forward problem, we prove that, under the natural fractional maximal-regularity condition $u\in W^{\alpha,p}(0,T;L^2(\Omega))\cap L^p(0,T;H^2(\Omega)\cap H_0^1(\Omega))$, the piecewise linear temporal reconstructions of PINN-loss minimizers converge to the exact solution. This result avoids the higher-order regularity assumptions on the exact solution that are often imposed in existing PINN convergence analyses. The proof combines $\liminf$/$\limsup$ arguments with Aubin--Lions compactness lemma. We also verify the required discrete maximal regularity and consistency properties for the aforementioned popular time discretizations. We then extend the analysis to the backward subdiffusion problem and establish convergence of the regularized recovery of the initial data from terminal observations, despite its ill-posedness. Numerical experiments for both forward and backward problems support the theoretical results, demonstrate the necessity of regularization, and show that $L^2$-type losses alone are generally insufficient. These results establish a rigorous theoretical basis for discretization-based PINNs for time-fractional diffusion problems.
\end{abstract}



\section{Introduction}
The subdiffusion model provides a fundamental modeling tool  for describing anomalous diffusion phenomena in complex systems, e.g., transport in porous media and biological processes. %
A wide variety of numerical methods have been developed for such problems (see, e.g., \cite{wu2025numerical,jin2025regularity,wang2020high,zhang2020numerical,li2023exponential,li2022exponential} for an incomplete list). Recently, deep learning-based solvers have emerged as a promising mesh-free alternative,  among which physics-informed neural networks (PINNs) \cite{raissi2019physics} have attracted considerable attention. PINNs approximate the solution of the PDE by a neural network (NN) trained through the minimization of a loss incorporating the residuals for the PDE, and initial and boundary conditions. By embedding the governing physical laws directly into the training process, PINNs avoid mesh generation and offer enhanced flexibility, especially for high-dimensional and inverse problems. Despite their empirical success, however, the theoretical understanding of PINNs remains incomplete, and rigorous convergence analysis is still an active area of research. 

Recent works have developed a unified framework for PINN error analysis based on stability, generalization, and optimization errors \cite{de2024numerical}. Within the framework, generalization bounds have been derived for PDEs and inverse problems \cite{mishra2023estimates,mishra2022estimates}, dimension-independent estimates have been obtained for several high-dimensional problems \cite{de2022error,de2024error}, and approximation errors have been related to the training loss and quadrature error for second-order hyperbolic PDEs \cite{qian2023physics}. Convergence results have also been derived for incompressible Navier--Stokes equations under suitable regularity assumptions \cite{biswas2022error,de2024errorNS}; see also \cite{zeinhofer2025unified} for an abstract framework. These works focus on deriving quantitative bounds under suitable \textit{a priori} regularity assumptions.

The second line of work focuses on the convergence of the PINN approximations.  Shin et al. \cite{shin2020convergence} proved strong $C^0$ convergence for second-order elliptic and parabolic equations via a Schauder-type argument and H\"older regularization. More recently, one may analyze the convergence in the framework of $\Gamma$-convergence, under weaker regularity assumptions. 
In this direction, \cite{gazoulis2025stability} established a stability and convergence framework for PINNs applied to elliptic and parabolic problems by combining coercivity, compactness, and $\Gamma$-convergence arguments, thereby proving convergence of minimizers to weak solutions in critical spaces. This variational perspective is closely connected to tools from nonlinear PDE analysis, e.g., weak compactness, lower semicontinuity, and energy estimates, and is particularly effective when classical high-order smoothness assumptions are unavailable. For time-dependent problems, maximal $L^p$ regularity provides a natural framework for controlling temporal derivatives and parabolic smoothing effects. By discretizing the time derivative via Runge--Kutta methods, Akrivis et al \cite{akrivis2025runge} proved the convergence of PINN minimizers for parabolic equations under discrete maximal regularity. This approach is especially relevant for nonlinear parabolic problems, where compactness-based techniques and maximal regularity often play a central role. However, the existing literature focuses predominantly on classical integer-order PDEs, and the theory for fractional evolution equations remains far less developed.

In this work, we develop a convergence analysis of PINNs for linear subdiffusion models. The fractional derivative involved in such models presents two inherent difficulties:  the lack of a useful chain rule and their nonlocal structure with weakly singular kernels. These features make the exact evaluation of the fractional derivative $\partial_t^\alpha u$ at training points computationally costly and often inaccurate. Inspired by \cite{akrivis2025runge,gazoulis2025stability}, and in the spirit of fPINNs \cite{pang2019fpinns}, we adopt a time-discrete treatment of the fractional derivative. We formulate an abstract convergence framework for time-discrete PINNs based on two key ingredients: discrete maximal regularity and consistency of the underlying time-stepping scheme. This framework applies to any time discretization satisfying these assumptions and naturally connects deep learning-based solvers with established tools from numerical analysis. As concrete examples, we consider several popular schemes, including backward Euler (BE), BDF2, L1, and fractional Crank--Nicolson (CN) schemes. Under the natural fractional maximal-regularity condition
$u\in W^{\alpha,p}(0,T;L^2(\Omega))\cap L^p(0,T;H^2(\Omega)\cap H_0^1(\Omega))$,
we rigorously prove the convergence of the resulting PINN schemes. This avoids the higher-order regularity assumptions on the exact solution that are often imposed in existing PINN convergence analyses. The proof relies crucially on the discrete maximal regularity of the time-stepping schemes. By combining $\liminf$ and $\limsup$ arguments with Aubin--Lions compactness lemma, we show that the piecewise linear temporal reconstructions of PINN-loss minimizers converge to the exact solution $u$ in $L^p(0,T;X)$. Finally, we verify the discrete maximal regularity and consistency assumptions required by the abstract framework for the aforementioned time discretizations.

We further extend the analysis to the backward subdiffusion problem, namely the recovery of the initial data $u_0$, and hence the reconstruction of the solution $u(t)$ for $0<t<T$, from the terminal observation $u(T)$.
This problem is motivated by scenarios in which one observes the current state of a subdiffusive process and seeks to reconstruct its earlier profile. This inverse problem has been studied from both analytical \cite{sakamoto2011initial,wang2013total} and numerical \cite{zhang2020numerical,zhang2023stability} perspectives. Beyond its intrinsic interest, the fractional backward problem may also be used as a regularization strategy for backward parabolic problems, as investigated in \cite{KaltenbacherRundell:2019}. Due to the ill-posed nature of the inverse problem, appropriate regularization must be incorporated into the PINN formulation. By establishing discrete regularity estimates for the backward problem, we prove convergence of the resulting regularized PINN approximation. Our analysis reveals that an $L^2$-type loss alone is generally insufficient to ensure convergence of the PINN approximation, and that additional regularization terms are needed, especially in the backward setting. This conclusion is supported by the numerical experiments in Section~\ref{numerics}.

The rest of the paper is organized as follows. In Section~\ref{pbm-stm}, we describe the subdiffusion problem, its temporal discretizations, and preliminary results from fractional calculus. In Section~\ref{pin}, we present the PINN formulation, develop the $\Gamma$-convergence framework, and prove the convergence for the forward problem. In Section~\ref{scm-asp}, we verify the required assumptions for several time-stepping schemes. In Section~\ref{cvg-bsp}, we extend the analysis to backward subdiffusion. Numerical experiments for both forward and backward problems are given in Section~\ref{numerics}. Finally, in Section~\ref{ccls}, we discuss possible future research directions. Throughout, let $Q_T=\Omega\times(0,T)$, and the notation $(\cdot,\cdot)_{L^2(D)}$, for $D\in\{(0,T), \Omega, Q_T\}$,  denotes the $L^2(D)$ inner product or suitable duality pairing. The notation $c$, with or without a subscript, denotes a generic positive constant which may change at each occurrence but is always independent of the neural network parameters and the time step size.

\section{Problem statement and preliminaries}\label{pbm-stm}
Let $\Omega\subset \mathbb{R}^d$ be a bounded domain with a smooth boundary $\partial\Omega$, and let $u_0 \in X \subset L^2(\Omega)$. Consider the following initial-boundary value problem of the linear subdiffusion model:
\begin{equation}\label{eq:subdiff}
\left\{\begin{aligned}
\partial_t^\alpha u - \Delta u &= f,\quad \mbox{in }\Omega\times(0,T],\\
u &= 0,\quad \mbox{on } \partial\Omega\times(0,T],\\
u(\cdot,0) &= u_0,\quad \mbox{in }\Omega,
\end{aligned}\right.
\end{equation}
where $u_0$ is the initial data, $f$ is the given smooth source, and $\partial_t^\alpha u$ denotes the Caputo fractional derivative in time of order $\alpha\in(0,1)$, defined by \cite[p. 42]{jin2021fractional}
\begin{equation}\label{deltaf-Caputo}
\partial_t^\alpha u(x,t) := \frac{1}{\Gamma(1-\alpha)}\int_0^t (t-s)^{-\alpha} \partial_s u(x,s)\,\mathrm{d}s.
\end{equation}

The natural maximal-regularity space for the solution $u$ to the model \eqref{eq:subdiff} is given by
\begin{align*}
    \mathfrak{S}:= L^p(0,T;H^{2}(\Omega)\cap H^{1}_0(\Omega)) 
    \cap W^{\alpha,p}(0,T;L^2(\Omega)).
\end{align*}
To ensure the maximal $L^p$ regularity of the solution $u$ to problem \eqref{eq:subdiff}, we assume that the initial data $ u_0 $ belongs to the space $ X = (L^{2}(\Omega),H^{2}(\Omega)\cap H^{1}_0 (\Omega))_{1-\frac{1}{\alpha p},p} $,
where $(\cdot,\cdot)_{\theta,p}$ denotes the real interpolation \cite{bergh1976interpolation}.

\subsection{Temporal discretization} 
We now introduce the time discretization for problem \eqref{eq:subdiff}. Let
$t_n=n\tau$, $n=0,1,\ldots,N$, be a uniform partition of $[0,T]$, with a time
step size $\tau=T/N$.
A time discretization scheme for approximately solving problem \eqref{eq:subdiff} reads: find $\tilde{u}^n\in H^2(\Omega)$, $n=1,2,\ldots,N$, such that 
\begin{equation}\label{eq:subdiff-discrete}
\left\{\begin{aligned}
 \bar{\partial}_\tau^\alpha u  - \Delta \tilde{u}^n &= \tilde{f}^n,\quad \text{in } \Omega\times(0,T],\\
\tilde{u}^n &= 0,\quad \text{on } \partial\Omega\times(0,T],\\
u^0 &= u_0,\quad \text{in } \Omega,
\end{aligned}\right.
\end{equation}
where $\bar{\partial}_\tau^\alpha u$ denotes a discrete approximation of the Caputo derivative $\partial_t^\alpha u(t)$:
\begin{align}\label{partial-tau-alpha}
    \bar{\partial}_\tau^\alpha u (t) = \tau^{-\alpha}\sum_{i=0}^n  a^{(n)}_{n-i} (u^i-u^0),\quad t \in (t_{n-1},t_n].
\end{align}  
The operator $\tilde{\cdot}$ denotes a suitable interpolation operator mapping the discrete sequence $\{u^n\}_{n=0}^N$ into $L^2(0,T;L^2(\Omega))$. %
The piecewise linear interpolant $\hat{u}(t)$ and the piecewise constant interpolant $\bar{u}(t)$ associated with a sequence $\{u^n\}_{n=0}^N$ are defined by
\begin{align*}
\hat{u}(t) = &\; \tfrac{t-t_{n-1}}{\tau} u^n + \tfrac{t_n - t}{\tau} u^{n-1},\quad t\in (t_{n-1},t_n], \\
\bar{u}(t) = &\; u^n,\quad t\in (t_{n-1},t_n],
\end{align*}
respectively. The coefficients $a^{(n)}_k$ and the interpolation operators $\tilde{u}$ for several commonly used schemes are listed below:
\begin{enumerate}
    \item BE scheme:  $a_k^{(n)} = a_k$ for $0 \le k \le n$, $ \sum_{k=0}^{\infty} a_k \xi^k = (1-\xi)^\alpha $, and $ \tilde{u} = \bar{u} $.
    \item BDF2 scheme: b $a_k^{(n)} = a_k$ for $0 \le k \le n$, $ \sum_{k=0}^{\infty} a_k \xi^k = (\frac{3}{2} - 2 \xi + \frac{1}{2} \xi^2)^\alpha $, and $ \tilde{u} = \bar{u} $. 
    \item L1 scheme: $ a_0 = \frac{1}{\Gamma (2-\alpha)} $, $ a_k^{(n)} = \frac{1}{\Gamma(2-\alpha)}\big((k+1)^{1-\alpha} - 2 k^{1-\alpha} + (k-1)^{1-\alpha}\big) $ for $ 1 \le k \le n-1 $, $a^{(n)}_n = \frac{1}{\Gamma (2-\alpha)}\big( (n-1)^{1-\alpha} - n^{1-\alpha} \big)$, and $ \tilde{u} = \bar{u} $.
    \item Fractional CN scheme: $a^{(n)}_k$ is the same as the BE scheme, and $ \tilde{u} (t) = \frac{\alpha}{2} u^{n-1} + (1-\frac{\alpha}{2})u^n $ for $ t \in (t_{n-1},t_n] $.
\end{enumerate} 

We define a fractional reconstruction for a sequence $\{u^n\}_{n=0}^{N} \in L^{2}(\Omega)$ by
\begin{align}\label{recons_frac_R}
    \mathcal{R}_{\tau}^{\alpha} u (t) := u^0 + \frac{1}{\Gamma (\alpha)} \int_{0}^{t} (t-s)^{\alpha-1} \bar{\partial}_{\tau}^{\alpha} u (s)\d s.
\end{align} 

The convergence analysis requires that the discrete  operator $\bar{\partial}_\tau^\alpha$ must approximate the continuous one $\partial_t^\alpha$ consistently and satisfy discrete maximal regularity. Thus we impose the following assumptions, which hold for the above four schemes.
\begin{assumption}\label{assump0}
	Fix $N$, and suppose $\{\tilde{u}^n\}_{n=1}^N \subset H^{1}_0 (\Omega) \cap H^2(\Omega)$ is the interpolation of $ u \in W^{1,p} (0,T;H^{1}_0 (\Omega)\cap H^{2}(\Omega)) $ and $ u^0 \in X $. The interpolation operator $ \tilde{\cdot} $ is bounded from $ W^{1,p} (0,T;H^{1}_0 (\Omega)\cap H^{2}(\Omega)) $ to $ L^p (0,T;H^{1}_0  (\Omega)\cap H^{2}(\Omega)) $, and satisfies 
\begin{align}
	\tilde u\to u\text{ in }L^p (0,T;H^{1}_0  (\Omega)\cap H^{2}(\Omega)) \text{ as } \tau \to 0. \label{eq:assump-tilde}
\end{align} 
	The discrete Caputo derivative $\bar{\partial}_\tau^\alpha$ satisfies the discrete maximal regularity estimate: 
\begin{align}\label{eq:assump0}
\| \bar{\partial}_{\tau}^\alpha u \|_{L^{p}(0,T;L^{2}(\Omega))} + \| \Delta \tilde{u} \|_{L^{p}(0,T;L^{2}(\Omega))} \le c \, \big(\| \bar{\partial}_{\tau}^\alpha u  -  \Delta \tilde{u} \|_{L^{p}(0,T;L^{2}(\Omega))} + \| u^0 \|_X\big), 
\end{align} 
    for some $ p > \frac{1}{\alpha} $.
    Additionally, $ \bar{\partial}_{\tau}^{\alpha} $ and $ \mathcal{R}_\tau^\alpha $ satisfy the following approximation properties:
\begin{align}
		\| \bar{\partial}_{\tau}^\alpha u - \partial_{t}^{\alpha} \hat{u} \|_{L^{p}(0,T;L^{2}(\Omega))} & \le c \hat{\beta}_{\tau} \| \hat{u} \|_{W^{1,p} (0,T;L^{2}(\Omega))} ,\label{eq_assu_barpt}\\
		\| \tilde{u} - \mathcal{R}_{\tau}^{\alpha} u \|_{L^{p}(0,T;L^{2}(\Omega))} & \le c \beta_{\tau} \| \bar{\partial}_{\tau}^{\alpha} u \|_{L^{p} (0,T;L^{2}(\Omega))} , \label{eq:assump0-2} 
	\end{align}
with $ \hat{\beta}_{\tau}, \beta_{\tau} \to 0 $ as $ \tau \to 0 $.  Further, $ \hat{u} $ and $ \tilde{u} $ satisfy the following estimate:
\begin{align}\label{eq:ass3}
    \| \hat{u} - \tilde{u} \|_{L^{p}(0,T;X)} \le c \big(\| \tilde{u} - \tilde{u}(\cdot - \tau) \|_{L^{p}(\tau,T;X)} + \| \tilde{u} \|_{L^{p}(0,\tau;X)} + \tilde{\beta}_{\tau} \| u^0 \|_{X}\big),
\end{align} 
with $ \tilde{\beta}_{\tau} \to 0 $ as $ \tau \to 0 $.
\end{assumption}

\subsection{Preliminaries on fractional calculus} \label{plinry}
Now we present preliminary results from fractional calculus that will be used below.
The left-sided Riemann-Liouville fractional integral ${}_0I_t^\alpha u$ is defined by
\begin{align*}
    {}_0 I_t^\alpha u (t) = \frac{1}{\Gamma (\alpha)} \int_0^t (t-s)^{\alpha - 1} u (s) \, \d s . 
\end{align*}

Let $ W_{L}^{\alpha,p} (0,T) $ be the set of functions whose extension by zero to $ (-\infty,0) $ belongs to $ W^{\alpha,p} (-\infty,T) $. The following result gives the boundedness of the operator $ {}_0 I_t^\alpha $, which can be derived using the argument of \cite[Theorem 2.9]{jin2021fractional}.
\begin{lemma}\label{thm:boundedness-left-RL-integral}
$ {}_0 I_t^\alpha $ extends continuously from $ L^p (0,T;L^2(\Omega)) $ to $ W^{\alpha,p}_{L} (0,T;L^2(\Omega)) $.
\end{lemma} 

The left-sided Riemann-Liouville fractional derivative ${}_0^R\!D_t^\alpha u$ is defined by
\begin{equation}\label{eq:def-right-RL-1}
{}_0^R\!D_t^\alpha u = \frac{1}{\Gamma(1-\alpha)} \frac{\d}{\d t} \int_0^t (t-s)^{-\alpha} u(s)\, \d s.
\end{equation}
By the argument of \cite[Corollary 2.2]{jin2021fractional}, we have the following lemma.
\begin{lemma} \label{thm:boundedness-right-RL-derivative}
${}_0^R\!D_t^\alpha $ is continuous from $ W^{\alpha,p}_{L} (0,T;L^2(\Omega)) $ to $ L^p (0,T;L^2(\Omega)) $.
\end{lemma} 

By \cite[Theorem 2.12]{jin2021fractional}, the Caputo derivative $\partial_t^\alpha$ coincides with the Riemann-Liouville derivative ${}_0^R\!D_t^\alpha$ on $W^{\alpha,p}_{L}(0,T;L^2(\Omega))$, since
\begin{align*}
   \partial_t^\alpha u = {}_0^R\!D_t^\alpha (u - u_0) ,\quad \forall u \in AC ([0,T];L^2(\Omega)),
\end{align*} 
where $AC([0,T])$ denotes the space of absolutely continuous functions over $[0,T]$.
Thus we have the following corollary, that shows that $\| u \|_{W^{\alpha,p}(0,T;L^{2}(\Omega))}$ is equivalent to $\| \partial_t^\alpha u \|_{L^{p}(0,T;L^{2}(\Omega))}$ for functions $u \in W^{\alpha,p}_{L}(0,T;L^2(\Omega))$.

\begin{corollary}\label{corol-fractional-sobolev}
The Caputo fractional differential operator $ \partial_t^\alpha $ extends continuously from $ W^{\alpha,p}_{0,L} (0,T;L^2(\Omega)) $ to $ L^p (0,T;L^2(\Omega)) $, and \begin{align*}%
        \| u \|_{W^{\alpha,p} (0,T;L^{2}(\Omega))} \le c \| \partial_t^\alpha u \|_{L^{p}(0,T;L^{2}(\Omega))} ,\quad \forall u \in W^{\alpha,p}_{L} (0,T;L^2(\Omega)).    \end{align*} 
\end{corollary} 
\begin{proof}
    It suffices to consider $ u \in C_{L}^{\infty} (0,T;L^{2}(\Omega)) $ whose extension by zero to $ (-\infty,0) $ belongs to $ C^{\infty}(0,T;L^{2}(\Omega)) $. By \cite[Theorem 2.13]{jin2021fractional}, we have
    $u = {}_0 I_t^\alpha \partial_t^\alpha u$. 
    The desired estimate then follows from the boundedness of $ {}_0 I_t^\alpha $ in Lemma~\ref{thm:boundedness-left-RL-integral}.
\end{proof}

\subsection{Aubin-Lions lemma}
We recall the compactness criterion in $ L^p (0,T;B) $. 
\begin{lemma}\cite[Theorem 5]{simon1986compact}\label{lem:compactness-criterion}
Let $ T > 0 $ and $ 1 \le p \le \infty $. $ X \subset B \subset Y $ are Banach spaces. $ X \hookrightarrow B $ compactly and $ B \hookrightarrow Y $ continuously. A bounded subset $ F \subset L^{p} ([0,T];X) $ is relatively compact in $ L^{p} (0,T;B) $ $($and in $ C (0,T;B) $ if $ p = \infty )$ if the following condition is satisfied:
\begin{equation}\label{condition-1}
\lim_{h \to 0} \| \tau_h f - f \|_{L^{p} (0,T-h;Y)} = 0 , \quad \text{uniformly for } f \in F ,
\end{equation}
with $ (\tau_h f) (t) = f (t+h) $ for $ t \in (0,T-h) $.
\end{lemma}

\section{PINNs and convergence for the forward problem}\label{pin}
In this section, we formulate the PINN for approximating the solution $u$ to problem \eqref{eq:subdiff} and establish its convergence.

\subsection{PINNs} 
PINNs are deep learning-based methods for solving PDEs, where the solution $u$ is approximated by an NN $u_\theta$ and the PDE is enforced through the loss. Our results depend only on the approximation capacity of the function class. An  NN maps every point $ (x,t) \in Q_T $ to $ u_\theta (x,t) \in \mathbb{R} $, through a composition of affine transformations and nonlinear activation. More precisely, for a given depth $ L \in \mathbb{N} $, we define the NN by
\begin{align*}
    u_\theta (x) = &\; A_L \sigma (A_{L-1} \sigma ( \cdots \sigma (A_1 (x,t) + b_1) + b_{L-2}) + b_{L-1}) + b_L, %
\end{align*}
where $ A_i \in \mathbb{R}^{d_{i} \times d_{i-1}} $ and $ b_i \in \mathbb{R}^{d_i} $ are the weight matrices and bias vectors, respectively, and $ \sigma : \mathbb{R} \to \mathbb{R} $ is a nonlinear activation and applied to a vector componentwise. The parameters of the NN are denoted by $ \theta = \{A_i, b_i\}_{i=1}^L $. The total number of parameters is given by $ \sum_{i=1}^L (d_i d_{i-1} + d_i) $. The set of all NNs with the given structure $ \mathcal{N}(L,d_1 , \cdots , d_L) $ is denoted by $ V_{\mathcal{N}} $.

Boundary conditions may be imposed through a suitable ansatz for the NN. For example, to enforce the zero Dirichlet boundary condition, we can use a suitable distance function $ \Phi(x) $ \cite[Section 5.1.1]{sukumar2022exact}. Let
\begin{align*}
    V_{\mathcal{N},0} = \{ u_\theta (x,t) = \Phi(x) v_\theta (x,t) : v_\theta \in V_{\mathcal{N}} \} .
\end{align*} 
Every function in $ V_{\mathcal{N},0} $ automatically satisfies the zero Dirichlet boundary condition.

The continuous loss functional $ \mathcal{G}(v) $ is defined as follows for an appropriate function $ v $
\begin{align*}
    \mathcal{G}(v) := \| \partial_{t}^{\alpha} v - \Delta v - f \|_{L^{p}(0,T;L^{2}(\Omega))}^{p} + \| v (\cdot,0) - u_0 \|_{X}^{2}.
\end{align*} 

For problem \eqref{eq:subdiff}, consider the PINN based on the temporal discretization \eqref{eq:subdiff-discrete}. We define a residual functional $\mathcal{G}_\tau(v)$ by
\begin{equation}\label{PINN-scheme}
\mathcal{G}_{\tau}(v) = \| \bar{\partial}_{\tau}^\alpha v - \Delta \tilde{v} - \tilde{f} \|_{L^{p}(0,T;L^{2}(\Omega))}^{p} + \| v (\cdot,0) - u_0 \|_{X}^2.
\end{equation}
We obtain the PINN approximation $u_\theta$ by minimizing $\mathcal{G}_\tau(v)$ over the class $V_{\mathcal{N},0}$.

We take a sequence of NN spaces such that for each $\ell\in\mathbb{N}$ there is an associated space $V_{\mathcal{N},0}$, denoted by $V_\ell$, satisfying the following Sobolev density property: for each $w\in W^{1,p}(0,T;H^{1}_0 (\Omega) \cap H^{2} (\Omega))$, there exists $w_\ell\in V_\ell$ such that
\begin{align}
    \|w_\ell-w\|_{W^{1,p}(0,T;H^{1}_0 (\Omega)\cap H^{2}(\Omega))}\to0
    \quad\text{as }\ell\to\infty. \label{eq:nn-approximation-property}
\end{align}
Since the temporal interpolation operators used below are uniformly stable in these spaces, the approximants can be selected so that
\begin{align}
&\|\hat w_\ell-w\|_{W^{1,p}(0,T;L^2(\Omega))}
+\|\tilde w_\ell-w\|_{L^p(0,T;H^{1}_0 (\Omega)\cap H^2(\Omega))}
+\|w_\ell^0-w(0)\|_X\to0. \label{eq:nn-approximation-property-2}
\end{align}
Thus, \eqref{eq:nn-approximation-property}--\eqref{eq:nn-approximation-property-2} are qualitative density properties only; in particular, no prescribed approximation rate relative to a norm of $w$ is required.

With the space $ V_\ell $, we adopt the shorthand notation 
\begin{align*}
    \mathcal{G}_\ell (v_\ell) = \mathcal{G}_{\tau (\ell)} (v_\ell) ,\quad \forall v_\ell \in V_\ell , %
\end{align*}
where the time-step size $ \tau (\ell) \to 0 $ as $ \ell \to \infty $. We still denote $ \tau (\ell) $ by $ \tau $ below. Let $\{ u_\ell \}$ be a sequence of minimizers such that
\begin{align}
    u_\ell := \arg\min_{v_\ell \in V_\ell} \mathcal{G}_\ell (v_\ell)  . \label{eq:residual-convergence}
\end{align} 
We also assume that the approximation of the source is consistent, i.e.,
\begin{align*}
    \tilde f_\ell\to f\quad\text{in }L^p(0,T;L^2(\Omega)).
\end{align*}
The goal is to show that $\{ u_\ell\}_\ell $ converges to the solution $u$ of problem \eqref{eq:subdiff} as $ \ell \to \infty $.

\subsection{Maximal regularity, and the \texorpdfstring{$ \liminf $--$ \limsup $}{liminf--limsup} arguments}
The convergence analysis of the PINN for the subdiffusion problem relies on the maximal regularity of the time-discrete solution. We first establish the following result.
\begin{proposition}\label{prop-31}
For a sequence $\{v_\ell\in V_\ell\}_{\ell=1}^\infty$,  if $ \{\mathcal{G}_\ell (v_\ell) \}_\ell $ is uniformly bounded, then under Assumption~\ref{assump0}, $ \tilde{v}_{\ell} $ is relatively compact in $ L^{p}(0,T;X) $. 
\end{proposition} 
\begin{proof}
Since $ \mathcal{G}_\ell(v_\ell) $ is bounded, by the triangle inequality, we have
\begin{align*}
&\| \bar{\partial}_{\tau}^\alpha v_\ell - \Delta \tilde{v}_\ell \|_{L^{p}(0,T;L^{2}(\Omega))}^{p} + \| v_\ell (\cdot,0) \|_{X}^2\\
\le & c\mathcal{G}_\ell (v_\ell) + c\| \tilde{f} \|_{L^{p}(0,T;L^{2}(\Omega))}^{p} + c \| u_0 \|_{X}^2 \le  c.
\end{align*} 
By the discrete maximal regularity estimate \eqref{eq:assump0}, we obtain
\begin{equation}\label{MR-regularity-discrete-u}
\| \bar{\partial}_{\tau}^\alpha v_\ell \|_{L^{p}(0,T;L^{2}(\Omega))}^{p} + \| \Delta \tilde{v}_\ell \|_{L^{p}(0,T;L^{2}(\Omega))}^{p} \le c.
\end{equation}
Hence, by the elliptic regularity pickup, 
\begin{align*}
        \| \tilde{v}_\ell \|_{L^{p}(0,T;H^{2}(\Omega))}^{p}\le c \| \Delta \tilde{v}_\ell \|_{L^{p}(0,T;L^{2}(\Omega))}^{p} \le &c  . %
    \end{align*}
Next we verify condition \eqref{condition-1} for $ \tilde{v}_{\ell} $ with $ Y = L^{2}(\Omega) $.  
\begin{align*}
	\| \tau_h \tilde{v}_{\ell} - &\tilde{v}_{\ell} \|_{L^{p}(0,T-h;L^{2}(\Omega))} 
	\le  \| \tau_h (\tilde{v}_{\ell} - \mathcal{R}_{\tau}^{\alpha} v_{\ell}) \|_{L^{p}(0,T-h;L^{2}(\Omega))}\\
    &+ \| \tau_h \mathcal{R}_{\tau}^{\alpha}v_\ell - \mathcal{R}_{\tau}^{\alpha}v_{\ell} \|_{L^{p}(0,T-h;L^{2}(\Omega))} + \| \mathcal{R}_{\tau}^{\alpha}v_{\ell} - \tilde{v}_{\ell} \|_{L^{p}(0,T-h;L^{2}(\Omega))}.
\end{align*}
Since $ \partial_{t}^{\alpha} \mathcal{R}_{\tau}^{\alpha} v_{\ell} = \bar{\partial}_{\tau}^{\alpha} v_{\ell} $, by Corollary~\ref{corol-fractional-sobolev}, we have
\begin{align*}
	&\| \tau_h \mathcal{R}_{\tau}^{\alpha}v_\ell - \mathcal{R}_{\tau}^{\alpha}v_{\ell} \|_{L^{p}(0,T-h;L^{2}(\Omega))} \le c h^\alpha \| \mathcal{R}_{\tau}^{\alpha}v_{\ell} \|_{W^{\alpha,p}(0,T;L^{2}(\Omega))} \\
    \le& ch^\alpha \| \partial_{t}^{\alpha}\mathcal{R}_{\tau}^{\alpha}v_{\ell} \|_{L^{p}(0,T;L^{2}(\Omega))}= c h^\alpha \| \bar{\partial}_{\tau}^{\alpha} v_{\ell} \|_{L^{p}(0,T;L^{2}(\Omega))} \le c h^\alpha .
\end{align*} 
By condition \eqref{eq:assump0-2}, we have
\begin{align*}
    \| \tau_h \tilde{v}_{\ell} - \tilde{v}_{\ell} \|_{L^{p}(0,T-h;L^{2}(\Omega))} \le c (h^\alpha + \beta_{\tau}).
\end{align*} 
Since $ \beta_{\tau} \to 0 $ as $ \tau \to 0 $, for any $ \varepsilon > 0 $, there exists $ L \in \mathbb{N} $ such that when $ \ell > L $, $ c \beta_{\tau} < \frac{\varepsilon}{2} $. Since translations are continuous for fixed $ \tilde{v}_{\ell} $ in $ L^{p}(0,T;L^{2}(\Omega)) $, there exists $ h_1 $ such that 
\begin{align*}
    \| \tau_h \tilde{v}_{\ell} - \tilde{v}_{\ell} \|_{L^{p}(0,T-h;L^{2}(\Omega))} \le \varepsilon \quad \text{for } 1 \le \ell \le L, \text{and } 0 < h < h_1.
\end{align*} 
$ h_1 $ may be dependent on $ \varepsilon $, but independent on $ \ell $. Also there exists $ h_2 $ such that $ c h^\alpha < \frac{\varepsilon}{2} $ for all $ 0 < h < h_2 $. Hence, we have 
\begin{align*}
    \| \tau_h \tilde{v}_{\ell} - \tilde{v}_{\ell} \|_{L^{p}(0,T-h;L^{2}(\Omega))} \le \varepsilon \quad \text{for } 0 < h < \min \{h_{1}, h_2\},
\end{align*} 
which implies the condition \eqref{condition-1}. The conclusion then follows from Lemma~\ref{lem:compactness-criterion} and the fact $ H^{2} (\Omega)\cap H^{1}_0 (\Omega) $ embeds compactly into $ X $.
\end{proof} 

We next employ the $\liminf$--$\limsup$ argument to prove the convergence of the PINN for the direct subdiffusion problem. First we state the following $\liminf$ inequality. 
\begin{lemma}\label{lim-inf-lem}
Let $\{ \mathcal{G}_\ell (v_\ell) \}_\ell$ be uniformly bounded, and suppose that $ \tilde{v}_\ell \to v $ in $ L^{p}(0,T;X) $. Then under Assumption~\ref{assump0}, we have $ v \in \mathfrak{S} $ and
    \begin{align}
        \mathcal{G} (v) \le \liminf_{\ell \to \infty} \mathcal{G}_\ell (v_\ell) .
    \end{align} 
\end{lemma}
\begin{proof}
    By the weak lower semicontinuity of the norms, it suffices to prove  
    \begin{align}
        \bar{\partial}_\tau^\alpha v_\ell - \Delta \tilde{v}_\ell - \tilde{f} \rightharpoonup  \partial_t^\alpha v - \Delta v - f \quad & \text{in} \quad L^{p}(0,T;L^{2}(\Omega)) , \label{eq:lem-5-1} \\
        v_\ell^0 \rightharpoonup v (\cdot,0)  \quad & \text{in} \quad X . \label{eq:lem-5-1-b}
    \end{align}  
Since $ \partial_{t}^{\alpha} \mathcal{R}_{\tau}^{\alpha} v_{\ell} = \bar{\partial}_{\tau}^{\alpha} v_{\ell} $, it suffices to show that 
\begin{align*}
\partial_t^\alpha \mathcal{R}_{\tau}^{\alpha}v_\ell \rightharpoonup \partial_t^\alpha v \quad\mbox{and}\quad \Delta \tilde{v}_\ell \rightharpoonup \Delta v \quad\mbox{in}\quad  L^{p}(0,T;L^{2}(\Omega)) . 
\end{align*}
By Corollary~\ref{corol-fractional-sobolev} and \eqref{MR-regularity-discrete-u}, the sequence $\{ \mathcal{R}_{\tau}^{\alpha} v_{\ell} - v_\ell^0 \}$ is bounded in $W^{\alpha,p}_{0,L}(0,T;L^2(\Omega))$. Hence, there exist $v_1 \in W^{\alpha,p}_{0,L}(0,T;L^2(\Omega))$ and a subsequence, still denoted by $\mathcal{R}_{\tau}^{\alpha} v_{\ell}$, such that
     \begin{align}
        \mathcal{R}_{\tau}^{\alpha} v_{\ell} - v_\ell^0 \rightharpoonup v_1 \quad \text{in } W^{\alpha,p}_{0,L} (0,T;L^2 (\Omega)). \label{con-vat-v1}
    \end{align}
Moreover, the boundedness of $\mathcal{G}_\ell(v_\ell)$ implies the uniform boundedness of $v_\ell^0$ in $X$. Thus, there exists $v_0 \in X$ such that $v_\ell^0 \rightharpoonup v_0$ in $X$. Since $X \subset L^{2}(\Omega)$, we have $\mathcal{R}_{\tau}^{\alpha} v_{\ell} \rightharpoonup v_0 + v_1$ in $W^{\alpha,p}(0,T;L^2(\Omega))$. Thus, $\mathcal{R}_{\tau}^{\alpha} v_{\ell} \rightharpoonup v_0 + v_1$ in $L^p(0,T;L^{2}(\Omega))$. The inequality \eqref{eq:assump0-2} indicates $ \mathcal{R}_{\tau}^{\alpha} v_{\ell} \to v $ in $ L^{p}(0,T;L^{2}(\Omega)) $. By the uniqueness of the limit, we have $v = v_0 + v_1 \in W^{\alpha,p}(0,T;L^2(\Omega))$. Since $p > \frac{1}{\alpha}$, the trace $v(\cdot,0)$ of $v$ (at $t=0$) is well-defined and $v(\cdot,0) = v_0$. This proves \eqref{eq:lem-5-1-b}.
From \eqref{con-vat-v1}, we have $\partial_t^\alpha (\mathcal{R}_{\tau}^{\alpha}{v}_\ell - v_\ell^0) \rightharpoonup \partial_t^\alpha (v - v_0)$ in $L^{p}(0,T;L^{2}(\Omega))$ by Corollary~\ref{corol-fractional-sobolev}. Hence, $\partial_t^\alpha \mathcal{R}_{\tau}^{\alpha}{v}_\ell \rightharpoonup \partial_t^\alpha v$ in $L^{p}(0,T;L^{2}(\Omega))$.
The convergence $\Delta \tilde{v}_\ell \rightharpoonup \Delta v$ in $L^{p}(0,T;L^{2}(\Omega))$ follows similarly from the above argument.
\end{proof} 

Next, we state the $\limsup$ inequality.
\begin{lemma}\label{lim-sup-lem}
    Let $ w \in \mathfrak{S} $. Under Assumption~\ref{assump0}, there exists a sequence $ \{ w_\ell \} $ such that $ w_\ell \to w $ in $ \mathfrak{S} $, $ w_{\ell} (\cdot,0) \to w (\cdot,0) $ in $ X $, and
    \begin{align*}
        \mathcal{G} (w) = \lim_{\ell \to \infty} \mathcal{G}_\ell (w_\ell) .
    \end{align*} 
    The sequence $ \{ w_{\ell} \} $ is called a recovery sequence for $ w $.
\end{lemma} 
\begin{proof}
By standard temporal mollification, $W^{1,p}(0,T;H^{1}_0 (\Omega)\cap H^{2}(\Omega))$ is dense in $\mathfrak S$. Since the trace map from $\mathfrak S$ to $X$ is continuous, there are $w_k\in W^{1,p}(0,T;H^{1}_0 (\Omega)\cap H^{2}(\Omega))$ such that
\begin{align}
    \|w_k-w\|_{\mathfrak S}+\|w_k(0)-w(0)\|_X\to0 \quad \text{as } k \to \infty. \label{eq:recovery-smooth}
\end{align}
Consequently, $\mathcal G(w_k)\to\mathcal G(w)$.

Fix $k$. By \eqref{eq:nn-approximation-property}--\eqref{eq:nn-approximation-property-2}, there are $w_{k,\ell}\in V_\ell$ such that $w_{k,\ell}\to w_k$ in $W^{1,p}(0,T;H^{1}_0 (\Omega)\cap H^{2}(\Omega))$. The piecewise linear interpolation $ \hat{\cdot} $ is uniformly bounded on $W^{1,p}(0,T;H^{1}_{0}(\Omega)\cap H^{2}(\Omega))$ and converges strongly to the identity as $\tau\to0$. By \eqref{eq:assump-tilde}, we have 
\begin{align}
    \hat w_{k,\ell}&\to w_k &&\text{in }W^{1,p}(0,T;L^2(\Omega)), \notag\\
    \tilde w_{k,\ell}&\to w_k &&\text{in }L^p(0,T;H^{1}_0 (\Omega)\cap H^{2}(\Omega)), \label{eq:lem-6-1}\\
    w_{k,\ell}^0&\to w_k(0) &&\text{in }X. \notag
\end{align}
Here the last convergence follows from $W^{1,p}(0,T;H^{1}_0 (\Omega)\cap H^{2}(\Omega))\hookrightarrow C([0,T];H^{1}_0 (\Omega)\cap H^{2}(\Omega))$ and $H^{1}_0 (\Omega)\cap H^{2}(\Omega)\hookrightarrow X$.

By condition \eqref{eq_assu_barpt} and the boundedness of ${}_0I_t^{1-\alpha}$ on $L^p(0,T;L^2(\Omega))$, we obtain
\begin{align*}
&\|\bar\partial_\tau^\alpha w_{k,\ell}-\partial_t^\alpha w_k\|_{L^p(0,T;L^2(\Omega))}\\
&\quad\le c\hat\beta_\tau\|\hat w_{k,\ell}\|_{W^{1,p}(0,T;L^2(\Omega))}
 +c\|\hat w_{k,\ell}'-w_k'\|_{L^p(0,T;L^2(\Omega))}\to0.
\end{align*}
Moreover, \eqref{eq:lem-6-1} implies $\Delta\tilde w_{k,\ell}\to\Delta w_k$ in $L^p(0,T;L^2(\Omega))$. Together with the source consistency $\tilde f_\ell\to f$ and the convergence of the initial values, this gives
\begin{align}
    \mathcal G_\ell(w_{k,\ell})\to\mathcal G(w_k)
    \quad\text{for every fixed }k. \label{eq:recovery-fixed-k}
\end{align}

It remains to make a diagonal choice. Choose an increasing sequence $L_k$ such that, for every $\ell\ge L_k$,
\begin{align*}
&\|w_{k,\ell}-w_k\|_{\mathfrak S}
+\|w_{k,\ell}^0-w_k(0)\|_X
+|\mathcal G_\ell(w_{k,\ell})-\mathcal G(w_k)|\le k^{-1}.
\end{align*}
For $L_k\le\ell<L_{k+1}$, set $w_\ell:=w_{k,\ell}$. Then $k\to\infty$ as $\ell\to\infty$, and \eqref{eq:recovery-smooth}--\eqref{eq:recovery-fixed-k} imply
$w_\ell\to w$ in $\mathfrak S$, $w_\ell^0\to w(0)$ in $X$, and $\mathcal G_\ell(w_\ell)\to\mathcal G(w)$.
\end{proof} 

\subsection{Convergence for the forward problem}\label{cvg-ana-fp}
We now state the first main result: the convergence of the PINN for the forward problem.
\begin{theorem}\label{main-thm-forward}
Let $ \{u_\ell\}_\ell $ be the PINN solutions defined in \eqref{eq:residual-convergence}, and let Assumption~\ref{assump0} hold. Then the piecewise linear temporal reconstructions $\{ \hat{u}_\ell \}_\ell $  converge to the solution $ u $ of \eqref{eq:subdiff} in $ L^{p}(0,T;X) $.
\end{theorem}
\begin{proof}
Let $\{ v_\ell \}_\ell$ be a recovery sequence for the exact solution $ u $. By Lemma~\ref{lim-sup-lem},
    \begin{align*}
        \mathcal{G}_\ell (v_\ell) \to \mathcal{G} (u) = 0 .
    \end{align*}
Since $ u_\ell $ minimizes the functional $ \mathcal{G}_\ell $, it follows that
$\mathcal{G}_\ell (u_\ell) \le \mathcal{G}_\ell (v_\ell) $.
Thus $ \mathcal{G}_\ell (u_\ell) $ is uniformly bounded with respect to $ \ell $.
By Proposition~\ref{prop-31},  the sequence $ \{\tilde{u}_\ell\} $ is relatively compact in $ L^{p}(0,T;X) $, there exists a subsequence of $ \{ u_\ell \} $, still denoted by $ \{ u_\ell \} $, and some $ u^* \in L^{p}(0,T;X) $ such that
\begin{equation*}
\tilde{u}_\ell \to u^* \quad \text{in} \quad L^{p}(0,T;X) .
\end{equation*}
To show that $ u^* $ coincides with the solution $ u $ of problem \eqref{eq:subdiff}, we claim that $ u^* $ minimizes the functional $ \mathcal{G} $. Indeed, by Lemmas~\ref{lim-inf-lem} and \ref{lim-sup-lem}, for any $w\in\mathfrak{S}$ with the recovery sequence ${w_\ell}$, we have 
\begin{align*}
    \mathcal{G} (u^*) \le &\; \liminf_{\ell \to \infty} \mathcal{G}_\ell (u_\ell) \le \limsup_{\ell \to \infty} \mathcal{G}_\ell (u_\ell) \le \limsup_{\ell \to \infty} \mathcal{G}_\ell (w_\ell) = \mathcal{G} (w),
\end{align*} 
Thus, $ u^* $ is the solution of problem \eqref{eq:subdiff}, which is unique by its well-posedness. Since every subsequence contains a sub-subsequence that converges in $L^p(0,T;X)$ to the solution $u^*$, a standard subsequence argument implies that the whole sequence converges.
It remains to show that $ \hat{u}_{\ell} - \tilde{u}_{\ell} \to 0 $ in $ L^{p}(0,T;X) $. By condition \eqref{eq:ass3}, when $ \ell $ is sufficiently large, we have
\begin{align*}
	\| \hat{u}_{\ell} - \tilde{u}_{\ell} \|_{L^{p}(0,T;X)} & \le c \big(\|\tilde{u}_{\ell} - \tilde{u}_{\ell} (\cdot - \tau)\|_{L^{p}(\tau,T;X)} + \| \tilde{u}_{\ell} \|_{L^{p}(0,\tau;X)} + \tilde{\beta}_{\tau}\| u_{\ell}^{0} \|_{X}\big) \\
	&\le c \big(\|\tilde{u}_{\ell} - \tilde{u}_{\ell} (\cdot - \tau)\|_{L^{p}(\tau,T;X)} + \| u \|_{L^{p}(0,\tau;X)} + \tilde{\beta}_{\tau}\| u_0 \|_{X}\big).
\end{align*} 
The term $ \|\tilde{u}_{\ell} - \tilde{u}_{\ell} (\cdot - \tau)\|_{L^{p}(\tau,T;X)}$ tends to $ 0 $ as $ \ell \to \infty $ since $ \tilde{u}_{\ell} $ converges to $ u $ in $ L^{p}(0,T;X) $. Hence, the right-hand side of the inequality tends to $ 0 $ as $ \ell \to \infty $.
\end{proof}

\section{Time-stepping schemes satisfying Assumptions~\ref{assump0}}\label{scm-asp}
In this section, we verify that the BE, BDF2, L1, and fractional CN schemes satisfy Assumptions~\ref{assump0}. 

\subsection{BE and BDF2 schemes}\label{subsec-BE-BDF2}
For the BE and BDF2 schemes, the discrete maximal regularity estimate \eqref{eq:assump0} has already been established in \cite{jin2018discrete} and \eqref{eq:assump-tilde} follows from $ \tilde{u} = \bar{u} $. Hence it remains to verify conditions \eqref{eq_assu_barpt}, \eqref{eq:assump0-2} and \eqref{eq:ass3}.  
We define the generating function
\begin{align}
    \delta_{\tau}^k (\xi) := \frac{\delta^k (\xi)}{\tau} , \quad \text{with } \delta^k (\xi) = \sum_{j = 1}^k \frac{1}{j} (1-\xi)^j .
\end{align} 
For the BE scheme, $ k = 1 $, while for the BDF2 scheme, $ k = 2 $. The following lemma states approximation properties of the generating function $ \delta_{\tau}^k (\xi) $ \cite[Lemma 3.1]{jin2023numerical}.
\begin{lemma}\label{prop-generating-func}
    For any $ \varepsilon > 0 $, there exists $ \theta_\varepsilon \in (\frac{\pi}{2},\pi) $ such that for any $ \theta \in (\frac{\pi}{2}, \theta_\varepsilon) $, there exist $ c$, $c_1 $, and $ c_2 $ such that
    \begin{align*}
        c_1 |z| \le |\delta_{\tau}^k (e^{-z\tau})| \le c_2 |z| , \quad |\delta_{\tau}^k(e^{-z \tau}) - z| \le c \tau^k |z|^{k+1} ,  \quad \forall z\in \Gamma_{\theta,\sigma}^\tau,\\
        |\delta_{\tau}^k(e^{-z \tau})^\alpha - z^\alpha| \le c \tau^k |z|^{\alpha + k} , \quad \forall z \in \Gamma_{\theta,\sigma}^\tau , 
    \end{align*} 
    where $ \sigma > 0 $ and the contour $ \Gamma_{\theta,\sigma}^\tau \subset \mathbb{C} $ is defined by
    \begin{align*}
        \Gamma_{\theta,\sigma}^\tau =  \{ z = \rho e^{\pm i \theta} : \rho \ge \sigma, | \Im (z) | \le \tfrac{\pi}{\tau} \} \cup \{ z = \sigma e^{i \varphi} : |\varphi| \le \theta \} .
    \end{align*} 
\end{lemma} 
By Lemma \ref{prop-generating-func}, we show that the BE and BDF2 schemes satisfy condition \eqref{eq_assu_barpt}. 
\begin{proposition}\label{prop-52}
The BE and BDF2 schemes satisfy
\begin{align}
\| \bar{\partial}_{ \tau}^\alpha u - \partial_t^\alpha \hat{u} \|_{L^{p}(0,T;L^{2}(\Omega))} \le c \tau^{1-\alpha} |\log \tau| \| \hat{u} \|_{W^{1,p}(0,T;L^{2}(\Omega))} . \label{eq:prop-7-1}
\end{align} 
\end{proposition} 
\begin{proof}
 Let  
\begin{align*}
        u (\xi) = \sum_{n=0}^\infty u^n \xi^n , \quad \delta_{\tau} u^n = \frac{u^n - u^{n-1}}{\tau}, \quad (\delta_{\tau}u) (\xi) = \sum_{n=1}^\infty (\delta_{\tau} u^n) \xi^n .
    \end{align*}
Further, we also define the truncation projection $ P_n $ by
$P_n u (\xi) = \sum_{m=0}^n u^m \xi^m $. Note that these sums are finite, since we may assume $ u^n = 0 $ for $ n > N $. By the definition of $ \bar{\partial}_{\tau}^\alpha $, we have
\begin{align}    \sum_{n=0}^{\infty}\bar{\partial}_{\tau}^\alpha u^n \xi^n =  \delta_{\tau}^k ( \xi )^\alpha \big(u (\xi) - \frac{u^0}{1-\xi}\big) . \label{eq:prop-7-2}
\end{align}
Rewriting $ \delta_\tau u (\xi) $ as
\begin{align*}
\delta_{\tau} u (\xi) = \frac{1}{\tau}\sum_{n=1}^{\infty} u^n \xi^n - \frac{1}{\tau} \sum_{n=0}^{\infty} u^n \xi^{n+1} = \frac{1-\xi}{\tau} u (\xi) = \delta_{\tau}^1 (\xi) u (\xi) - \frac{u^0}{\tau},
\end{align*}   
it follows that 
\begin{align*}
\sum_{n=0}^{\infty} \bar{\partial}_{\tau}^\alpha u (t_n) \xi^n = \delta_{\tau}^{k} (\xi)^{\alpha} \delta_{\tau}^1 (\xi)^{-1} \delta_{\tau} u (\xi) . %
\end{align*} 
By Cauchy's integral formula and changing variables $\xi=e^{-z\tau}$, we get 
    \begin{align}
        \bar{\partial}_{\tau}^\alpha u (t_n) = &\; \frac{1}{2\pi i} \int_{|\xi| = \rho} \xi^{-n-1}\delta_{\tau}^{k} (\xi)^{\alpha} \delta_{\tau}^1 (\xi)^{-1} (P_n \delta_{\tau} u) (\xi)  \, \d \xi \notag \\
        = &\; \frac{\tau}{2\pi i} \int_{\Gamma^\tau}e^{z t_n}  \delta_{\tau}^{k} (e^{-z \tau})^{\alpha} \delta_{\tau}^1 (e^{-z \tau})^{-1} (P_n\delta_{\tau} u) (e^{-z \tau})  \, \d z , \label{eq:prop-7-4}
    \end{align}
with the contour $ \Gamma^\tau $ defined by 
$\Gamma^\tau = \{ z = -\frac{1}{\tau} \log \rho + i y : y \in [-\frac{\pi}{\tau},\frac{\pi}{\tau}] \} $,
    oriented with a decreasing imaginary part. By the analyticity of the generating function and periodicity of $ e^{-z \tau} $,  Cauchy's theorem allows deforming the contour $ \Gamma^{\tau} $ to $ \Gamma_{\theta,\sigma}^\tau $ oriented with an increasing imaginary part. Consequently,
    \begin{align}\label{exp-barpartial-tau}
        \bar{\partial}_{\tau}^\alpha u (t_n) = \frac{\tau}{2\pi i} \int_{\Gamma_{\theta,\sigma}^\tau}e^{z t_n}  \delta_{\tau}^{k} (e^{-z \tau})^{\alpha} \delta_{\tau}^1 (e^{-z \tau})^{-1} (P_n \delta_{\tau} u) (e^{-z \tau})  \, \d z .
    \end{align} 
Meanwhile, since the Laplace transform of $\frac{t^{-\alpha}}{\Gamma(1-\alpha)}$ is $ z^{\alpha-1}$, by taking the inverse Laplace transform, we have
\begin{align}\label{eq:tmal}
\frac{t^{-\alpha}}{\Gamma (1-\alpha)} & = \frac{1}{2 \pi i} \int_{\Re (z) = c} e^{z t} z^{\alpha - 1} \d z 
 = \frac{1}{2\pi i} \int_{\Gamma_{\theta,\sigma}} e^{z t} z^{\alpha - 1} \d z,
    \end{align}
with the contour $   \Gamma_{\theta,\sigma} =  \{ z = \rho e^{\pm i \theta} : \rho \ge \sigma \} \cup \{ z = \sigma e^{i \varphi} : |\varphi| \le \theta \} $, and $ c > 0 $.
Then by the definition of $ \partial_t^\alpha \hat{u} $, we have 
\begin{align}\label{alter-def-ptalu}
        \partial_t^\alpha \hat{u} (t_n) = \frac{1}{2\pi i}\int_0^{t_n} \hat{u}' (s) \int_{\Gamma_{\theta,\sigma}} e^{z (t_n - s)} z^{\alpha - 1} \d z \d s.
    \end{align}
Then by Fubini's theorem, we obtain 
    \begin{align}
        \partial_t^\alpha \hat{u} (t_n) = &\; \frac{1}{2\pi i} \int_{\Gamma_{\theta,\sigma}} e^{z t_n} z^{\alpha - 1} \int_{0}^{t_n} e^{-zs} \hat{u}' (s)   \, \d s \d z %
        \notag \\
        = &\; \frac{1}{2 \pi i} \int_{\Gamma_{\theta,\sigma}} e^{z t_n} z^{\alpha - 1} \sum_{j=1}^n \int_{t_{j-1}}^{t_j} e^{-z s} \delta_{\tau}u^j   \, \d s \d z \notag \\
        = &\; \frac{1}{2 \pi i} \int_{\Gamma_{\theta,\sigma}} e^{z t_n} z^{\alpha - 2} \sum_{j = 1}^n \delta_{\tau}u^j (e^{-z t_{j-1}} - e^{-z t_j})\, \d z \notag \\
        = &\; \frac{1}{2 \pi i} \int_{\Gamma_{\theta,\sigma}} e^{z t_n} z^{\alpha - 2} (e^{z \tau} - 1) (P_n\delta_{\tau}u) (e^{-z \tau}) \, \d z , \label{eq:prop-7-5}
    \end{align}
    Hence, for $ t \in (t_{n-1},t_n] $, we can decompose the error as
\begin{align*}
\bar{\partial}_{\tau}^\alpha u (t_n) - \partial_t^\alpha \hat{u} (t) = &\big(\bar{\partial}_{\tau}^\alpha u (t_n) - \partial_t^\alpha \hat{u} (t_n) \big)+\big( \partial_t^\alpha \hat{u} (t_n) - \partial_t^\alpha \hat{u} (t)\big) 
=: \mathrm{I}_1 (t) + \mathrm{I}_2 (t) . \end{align*} 
Combining \eqref{exp-barpartial-tau} and \eqref{eq:prop-7-5} gives the following splitting 
\begin{align}
        \mathrm{I}_1 (t) = &\; \frac{1}{2\pi i} \int_{\Gamma_{\theta,\sigma}^\tau} e^{z t_n} \big(\tau \delta_{\tau}^{k} (e^{-z \tau})^{\alpha} \delta_{\tau}^1 (e^{-z \tau})^{-1} - z^{\alpha - 2} (e^{z \tau}-1) \big) (P_n \delta_{\tau} u) (e^{-z \tau})  \, \d z \notag \\
        &\; - \frac{1}{2 \pi i} \int_{\Gamma_{\theta,\sigma}\setminus \Gamma_{\theta,\sigma}^\tau}\!\!\!\!\! e^{zt_n} z^{\alpha - 2} (e^{z \tau} - 1) (P_n\delta_{\tau}u) (e^{-z \tau}) \, \d z 
        =:  \mathrm{I}_{1,1} (t) + \mathrm{I}_{1,2} (t).\label{eq:prop-7-7}
    \end{align}
    By Lemma~\ref{prop-generating-func}, we have the estimate
    \begin{align*}
        |\tau\delta_{\tau}^{k} (e^{-z \tau})^{\alpha} \delta_{\tau}^1 (e^{-z \tau})^{-1} - z^{\alpha - 2} (e^{z \tau}-1)| \le c \tau^{2-\alpha }  , \quad \forall z \in \Gamma_{\theta,\sigma}^\tau . %
    \end{align*} 
This estimate implies
\begin{align}
|\mathrm{I}_{1,1} (t)|  \le c \tau^{2-\alpha } \sum_{j=1}^n |\delta_{\tau}u^j|\int_{\Gamma_{\theta,\sigma}^\tau} e^{\Re(z) t_{n-j}}  \, |\d z| . \label{eq:prop-7-9}
\end{align} 
By choosing $ \sigma = t_{n-j+1}^{-1} $ in the contour integral, we deduce
\begin{align} \int_{\Gamma_{\theta,\sigma}^\tau} \!\!\!\! e^{\Re(z) t_{n-j}}  \, |\d z| \le &c \int_{\sigma}^{\frac{\pi}{\tau \sin(\theta)}}\!\!\!\! e^{r t_{n-j} \cos(\theta)}  \, \d r + c \int_{-\theta}^{\theta} \!\!\!\!e^{\sigma t_{n-j} \cos(\varphi)} \sigma \, \d \varphi 
        \le  ct_{n-j+1}^{ - 1} . \label{eq:prop-7-10}
    \end{align} 
Substituting this estimate into \eqref{eq:prop-7-9} gives
\begin{align*}
|\mathrm{I}_{1,1} (t)| \le c \tau^{2-\alpha } \sum_{j=1}^n t_{n-j+1}^{ - 1} |\delta_{\tau} u^j| . %
\end{align*}
Similarly, for the term $ \mathrm{I}_{1,2} (t) $, direct computation yields
\begin{align*}
|\mathrm{I}_{1,2} (t)| \le &c \sum_{j=1}^n |\delta_{\tau} u^j| \int_{\frac{\pi}{\tau \sin{\theta}}}^{\infty} e^{t_{n-j}\rho \cos{\theta}} \rho^{\alpha-2} \, \d \rho \notag \\
        \le &c \tau^{2-\alpha}\sum_{j=1}^n |\delta_{\tau} u^j| \int_{\frac{\pi}{\tau \sin{\theta}}}^{\infty} e^{t_{n-j} \rho \cos{\theta}} \, \d \rho \notag \\
        \le &c \tau^{2-\alpha} \sum_{j=1}^n t_{n-j+1}^{-1} |\delta_{\tau}u^j|.
\end{align*} 
By combining the last two estimates, and using H\"{o}lder's inequality, we deduce
    \begin{align}
        \|\mathrm{I}_1 \|_{L^{p}(0,T;L^{2}(\Omega))}^{p} \le &  c \tau^{1 + p (2-\alpha)} \sum_{n=1}^N \Big| \sum_{j=1}^n t_{n-j+1}^{-1} \|\delta_{\tau} u^j \|_{L^{2}(\Omega)} \Big|^p \notag \\
        \le &c \tau^{1 + p (2-\alpha)} \sum_{n=1}^N \Big(\sum_{j=1}^n t_{n-j+1}^{-1}\Big)^{p-1} \Big( \sum_{j=1}^n t_{n-j+1}^{-1} \|\delta_{\tau}u^j \|_{L^{2}(\Omega)}^p \Big) \notag \\
        \le &c \tau^{2 + p (1-\alpha) } |\log \tau|^{p-1} \sum_{n=1}^N \sum_{j=1}^n t_{n-j+1}^{-1} \|\delta_{\tau}u^j\|_{L^{2}(\Omega)}^p \notag \\
        \le &c \tau^{1 + p (1-\alpha)} |\log \tau |^p \sum_{n=1}^N \|\delta_{\tau}u^n\|_{L^{2}(\Omega)}^{p} \notag \\
        = &c \tau^{p (1-\alpha)} |\log \tau|^p \| \hat{u}' \|_{L^{p}(0,T;L^{2}(\Omega))}^{p}. \label{eq:prop-7-13}
    \end{align}
    For the term $ \mathrm{I}_2 (t) $, by the definition \eqref{deltaf-Caputo} of the fractional derivative $ \partial_t^\alpha \hat{u} $, we have
    \begin{align*}
        \mathrm{I}_2 (t) = &\;  \frac{1}{\Gamma (1-\alpha)} \int_{0}^{t_n} (t_n-s)^{-\alpha} \hat{u}' (s) \d s - \frac{1}{\Gamma (1-\alpha)} \int_{0}^{t} (t-s)^{-\alpha} \hat{u}' (s) \d s \notag \\
        = &\; \frac{1}{\Gamma (1-\alpha)} \int_{t}^{t_n} (t_n-s)^{-\alpha} \hat{u}' (s) \d s + \frac{1}{\Gamma (1-\alpha)}\int_{0}^t \zeta (t,s) \hat{u}' (s) \d s \notag \\
        =: &\; \mathrm{I}_{2,1} (t) + \mathrm{I}_{2,2} (t). %
    \end{align*}
    with $ \zeta (t,s) := (t_n - s)^{-\alpha} - (t-s)^{-\alpha} $ for $ t \in (t_{n-1},t_n] $, $ s \in (0,t) $, which satisfies \begin{align}\label{eti-zeta}
        |\zeta (t,s)| \le \min \{\alpha \tau (t-s)^{-\alpha-1}, (t-s)^{-\alpha}\}.
    \end{align} 
    For the term $ \mathrm{I}_{2,1} (t) $, we obtain
    \begin{align}
        \| \mathrm{I}_{2,1} \|_{L^{p}(0,T;L^{2}(\Omega))}^{p} = &\; \frac{1}{\Gamma (1-\alpha)^p} \sum_{n=1}^N \int_{t_{n-1}}^{t_n} \Big( \int_{t}^{t_n} (t_n - s)^{-\alpha}   \| \delta_{\tau} u^n \|_{L^{2}(\Omega)} \, \d s\Big)^p\, \d t \notag \\
        = &\; \frac{1}{\Gamma (1-\alpha)^p} \sum_{n=1}^N \int_{t_{n-1}}^{t_n} \Big( \frac{(t_n - t)^{1-\alpha}}{1-\alpha}  \| \delta_{\tau} u^n \|_{L^{2}(\Omega)} \Big)^p \, \d t \notag \\
        \le &\; c \tau^{p (1 - \alpha)} \| \hat{u}' \|_{L^{p}(0,T;L^{2}(\Omega))}^{p}. \label{eq:prop-7-15}
    \end{align} 
For the term $ \mathrm{I}_{2,2} (t) $, by H\"{o}lder's inequality, we derive
    \begin{align*}
        &\| \mathrm{I}_{2,2} \|_{L^{p}(0,T;L^{2}(\Omega))}^{p} \le \; c \int_{0}^{T} \Big( \int_{0}^t |\zeta (t,s)| \| \hat{u}'(s) \|_{L^{2}(\Omega)}\, \d s   \Big)^p\, \d t \notag \\
        \le &\; c  \int_{0}^{T} \Big(\int_{0}^t |\zeta (t,s)|\,\d s \Big)^{p-1} \int_{0}^t |\zeta (t,s)| \| \hat{u}' (s) \|_{L^{2}(\Omega)}^p \d s \d t .
    \end{align*} 
    Using the estimate \eqref{eti-zeta}, we have
    \begin{align*}
		\int_{0}^t |\zeta (t,s)| \,\d s \le &\; \int_{0}^{\max{0,t-\tau}} \alpha \tau (t-s)^{-\alpha - 1} \, \d s + \int_{\max{0,t-\tau}}^{t} (t-s)^{-\alpha} \, \d s 
        \le c \tau^{1-\alpha} .
    \end{align*} 
Therefore, we arrive at
\begin{align*}
    \| \mathrm{I}_{2,2} \|_{L^{p}(0,T;L^{2}(\Omega))}^{p} \le c \tau^{(1-\alpha)(p-1)}  \int_{0}^{T} \int_{0}^t |\zeta (t,s)| \| \hat{u}' (s) \|_{L^{2}(\Omega)}^p \d s \d t .
\end{align*}
Changing the integration order  and applying the bound \eqref{eti-zeta} again yield
    \begin{align}
        \| \mathrm{I}_{2,2} \|_{L^{p}(0,T;L^{2}(\Omega))}^{p} \le &\; c \tau^{(1-\alpha)(p-1)} \int_{0}^T \| \hat{u}' (s) \|_{L^{2}(\Omega)}^p \int_{s}^{T} |\zeta (t,s)| \, \d t \d s \notag \\
        \le &\; c \tau^{(1-\alpha)p} \| \hat{u}' \|_{L^{p}(0,T;L^{2}(\Omega))}^{p} . \label{I_22-eti}
    \end{align}
  By combining the estimates \eqref{eq:prop-7-15}, and \eqref{I_22-eti}, we conclude
    \begin{align}
        \|\mathrm{I}_2 \|_{L^{p}(0,T;L^{2}(\Omega))} \le c \tau^{1-\alpha}  \| \hat{u}' \|_{L^{p}(0,T;L^{2}(\Omega))}. \label{eq:prop-7-18}
    \end{align}
    Finally, this estimate and \eqref{eq:prop-7-13} complete the proof.
\end{proof} 

Now we show that the BE and BDF2 schemes satisfy condition \eqref{eq:assump0-2} with $ \beta_{\tau} = \tau^\alpha |\log \tau| $. 
\begin{proposition}\label{prop-52-R}
The BE and BDF2 schemes satisfy
\begin{align*}
	\| \bar{u} - \mathcal{R}_{\tau}^{\alpha}u \|_{L^{p}(0,T;L^{2}(\Omega))} \le c \tau^{\alpha} | \log \tau | \| \bar{\partial}_{\tau}^{\alpha} u \|_{L^{p}(0,T;L^{2}(\Omega))} . %
\end{align*} 
\end{proposition} 
\begin{proof}
Let $ \bar{\mathcal{R}}_{\tau}^{\alpha} u (t) := \mathcal{R}_{\tau}^{\alpha} u (t_n) $ for $ t \in (t_{n-1},t_n ] $. Then by the interpolation error, we have 
\begin{align}
	\| \bar{\mathcal{R}}_{\tau}^{\alpha} u - \mathcal{R}_{\tau}^{\alpha}u \|_{L^{p}(0,T;L^{2}(\Omega))} \le c \tau^{\alpha} \| \partial_{t}^{\alpha}\mathcal{R}_{\tau}^{\alpha}u \|_{L^{p}(0,T;L^{2}(\Omega))} = c \tau^{\alpha} \| \bar{\partial}_{\tau}^{\alpha}u \|_{L^{p}(0,T;L^{2}(\Omega))}, \label{est_int_Ru_R}
\end{align} 
for $ \alpha p >1 $. Hence, it suffices to bound the error $ \bar{u} - \bar{\mathcal{R}}_{\tau}^{\alpha}u $. Similarly to \eqref{exp-barpartial-tau}, we get
    \begin{align}\label{exp-barpartial-tau-R}
		 u^n = u^0
		+\frac{\tau}{2\pi i} \int_{\Gamma_{\theta,\sigma}^\tau}e^{z t_n}
		\delta_{\tau}^{k} (e^{-z \tau})^{-\alpha}
		\sum_{j=1}^n \bar{\partial}_{\tau}^{\alpha} u^j e^{-z t_j} \, \d z .
    \end{align}
    Similarly to \eqref{eq:prop-7-5}, we get
        \begin{align}
        \mathcal{R}_{\tau}^{\alpha}u (t_n) = 
        u^0 + \frac{1}{2 \pi i} \int_{\Gamma_{\theta,\sigma}} e^{z t_n} z^{-\alpha - 1} (e^{z \tau} - 1) \sum_{j=1}^n \bar{\partial}_{\tau}^{\alpha}u^j e^{-z t_j} \, \d z , \label{eq:prop-7-5-R}
    \end{align}
    By Lemma~\ref{prop-generating-func}, we have the estimate
    \begin{align}
        |\tau\delta_{\tau}^{k} (e^{-z \tau})^{-\alpha}   - z^{-\alpha - 1} (e^{z \tau}-1)| \le c \tau^{1+\alpha }  , \quad \forall z \in \Gamma_{\theta,\sigma}^\tau . \label{eq:prop-7-8-Rerror}
    \end{align} 
By combining \eqref{exp-barpartial-tau-R}, \eqref{eq:prop-7-5-R}, and \eqref{eq:prop-7-8-Rerror}, using a similar splitting to \eqref{eq:prop-7-7}, and following the subsequent proof, we can derive the desired result.
\end{proof} 

Next we verify condition \eqref{eq:ass3} for BE and BDF2 schemes with $ \tilde{\beta}_{\tau} = \tau^{1/p} $.
\begin{proposition}\label{prop-55}
The interpolations $ \hat{u} $ and $ \bar{u} $ for a sequence $ \{u^n\}_{n=0}^{N} \subset X  $ satisfy
\begin{align}
    \| \hat{u} - \bar{u} \|_{L^{p}(0,T;X)} \le c \big(\| \bar{u} - \bar{u}(\cdot - \tau) \|_{L^{p}(\tau,T;X)} + \| \bar{u} \|_{L^{p}(0,\tau;X)} + \tau^{1/p} \| u^0 \|_{X}\big). \label{eq:prop-7-barhat}
\end{align} 
\end{proposition}
\begin{proof}
For $t\in(t_{n-1},t_n]$, the definitions of $\hat u$ and $\bar u$ give
\begin{align*}
\hat u(t)-\bar u(t)=\frac{t_n-t}{\tau}(u^{n-1}-u^n).
\end{align*}
Consequently,
\begin{align*}
\|\hat u-\bar u\|_{L^p(0,T;X)}^p
=\frac{\tau}{p+1}\sum_{n=1}^N\|u^n-u^{n-1}\|_X^p.
\end{align*}
For $n\ge2$, we have
\begin{align*}
\tau\sum_{n=2}^N\|u^n-u^{n-1}\|_X^p
=\|\bar u-\bar u(\cdot-\tau)\|_{L^p(\tau,T;X)}^p.
\end{align*}
For the first interval, the triangle inequality yields
\begin{align*}
\tau^{1/p}\|u^1-u^0\|_X
&\le \tau^{1/p}\|u^1\|_X+\tau^{1/p}\|u^0\|_X\\
&=\|\bar u\|_{L^p(0,\tau;X)}+\tau^{1/p}\|u^0\|_X.
\end{align*}
Combining the preceding three identities proves \eqref{eq:prop-7-barhat}.
\end{proof}

\subsection{L1 Scheme} 
For the L1 scheme, the discrete maximal regularity \eqref{eq:assump0} has been proved in \cite{jin2018discrete} and \eqref{eq:assump-tilde} follows from $ \tilde{u} = \bar{u} $. Since for the L1 scheme, $ \bar{\partial}_{\tau}^\alpha u (t) = \partial_t^\alpha \hat{u} (t_n) $, the condition \eqref{eq_assu_barpt} has been proved in \eqref{eq:prop-7-18}. Condition \eqref{eq:ass3} is the same as Proposition~\ref{prop-55}. Therefore, it suffices to prove condition \eqref{eq:assump0-2}. We define the function
\begin{align}\label{def_gamma_xi}
    \gamma_{\tau}(\xi)^\alpha := \frac{(1-\xi)^{2}\operatorname{Li}_{\alpha-1}(\xi)}{\tau^{\alpha}\xi \Gamma (2-\alpha)},
\end{align} 
where the polylogarithm function $ \operatorname{Li}_{p}(\xi) $ is defined by $ \operatorname{Li}_{p}(\xi):= \sum_{j=1}^{\infty}\frac{\xi^j}{j^p}. $ 
The following lemma states approximation properties of the function $ \gamma_{\tau} (\xi) $ \cite[Lemma 5.3, Corollary 5.1]{jin2023numerical}.
\begin{lemma}\label{prop-generating-func-L1}
    Let $\theta$ be close to $\frac{\pi}{2}$ and $\sigma < \frac{\pi}{2}\tau$. Then there exist $ c$, $c_1 $, and $ c_2 $ such that
    \begin{align*}
        c_1 |z| \le |\gamma_{\tau} (e^{-z\tau})| \le c_2 |z| , \quad 
        |\gamma_{\tau}(e^{-z \tau})^\alpha - z^\alpha| \le c \tau^{2-\alpha} |z|^{2} , \quad \forall z \in \Gamma_{\theta,\sigma}^\tau , 
    \end{align*} 
\end{lemma} 

\begin{proposition}\label{prop-54}
The L1 scheme satisfies
\begin{align*} 
	\| \bar{u} - \mathcal{R}_{\tau}^{\alpha}u \|_{L^{p}(0,T;L^{2}(\Omega))} \le c \tau^{\alpha} |\log \tau |  \| \bar{\partial}_{\tau}^{\alpha} u \|_{L^{p}(0,T;L^{2}(\Omega))} .
\end{align*} 
\end{proposition}
\begin{proof}
Using \eqref{est_int_Ru_R}, it suffices to estimate the error $ \bar{u} - \bar{\mathcal{R}}_{\tau}^{\alpha}u $. For the L1 scheme, we have
    \begin{align*}
        \bar{\partial}_{\tau}^\alpha u^n = &\; \frac{1}{\Gamma (1-\alpha)} \sum_{j=1}^n \int_{t_{j-1}}^{t_j} (t_n - s)^{-\alpha} \delta_{\tau} u^j \, \d s 
        = \; \sum_{j=1}^{n} b_{n-j} \delta_{\tau} u^j , %
    \end{align*} 
with $b_j := \frac{\tau^{1-\alpha}}{\Gamma (2-\alpha)} ( (j+1)^{1-\alpha} - j^{1-\alpha} )$. By \cite[Proposition 5.1]{jin2023numerical}, $ u^{n} $ has the following representation 
\begin{align*}
u^n = u^0+\frac{\tau}{2\pi i} \int_{\Gamma_{\theta,\sigma}^\tau}e^{z t_n}\gamma_{\tau} (e^{-z \tau})^{-\alpha}\sum_{j=1}^n \bar{\partial}_{\tau}^{\alpha} u^j e^{-z t_j} \, \d z .
\end{align*}
Using Lemma~\ref{prop-generating-func-L1}, we derive that for some appropriate $\theta, \sigma$,
\begin{align*}
        |\tau\gamma_{\tau} (e^{-z \tau})^{-\alpha}   - z^{-\alpha - 1} (e^{z \tau}-1)| \le c \tau^{1+\alpha }  , \quad \forall z \in \Gamma_{\theta,\sigma}^\tau . %
\end{align*} 
The repeating the argument of Proposition~\ref{prop-52} completes the proof.
\end{proof} 

\subsection{Fractional CN Scheme} 
For the fractional CN scheme, the discrete maximal regularity \eqref{eq:assump0} has been proved in \cite{jin2018discrete} and \eqref{eq:assump-tilde} follows from $ \tilde{u}^n = \frac{\alpha}{2}u^{n-1} + (1-\frac{\alpha}{2})u^n $. Condition \eqref{eq_assu_barpt} is the same as BE scheme. Let the translation $ \bar{u}(\cdot-\tau) $ take the value $ u^0 $ on the first interval $ (0,\tau] $. Then $ \tilde{u} = \frac{\alpha}{2}\bar{u}(\cdot-\tau) + (1-\frac{\alpha}{2})\bar{u} $. By Proposition~\ref{prop-52} and a translation estimate for $ \mathcal{R}_{\tau}^{\alpha}u $, we can derive that the fractional CN scheme satisfies condition \eqref{eq:assump0-2}.

	Next, we verify condition \eqref{eq:ass3} for fractional scheme with $ \tilde{\beta}_{\tau} = \tau^{1/p} $.
	\begin{proposition}
The interpolations $ \hat{u} $ and $ \tilde{u} $ for a sequence $ \{u^n\}_{n=0}^{N} \subset X  $ satisfy
\begin{align}
    \| \hat{u} - \tilde{u} \|_{L^{p}(0,T;X)} \le c \big(\| \tilde{u} - \tilde{u}(\cdot - \tau) \|_{L^{p}(\tau,T;X)} + \| \tilde{u} \|_{L^{p}(0,\tau;X)} + \tau^{1/p} \| u^0 \|_{X}\big). 
	\end{align} 
		\end{proposition} 
\begin{proof}
For $t\in(t_{n-1},t_n]$, the definitions of $\hat u$ and $\tilde u$ give
\begin{align*}
\hat u(t)-\tilde u(t)
=\Big(\frac{t-t_{n-1}}{\tau}-1+\frac\alpha2\Big)
(u^n-u^{n-1}).
\end{align*}
Consequently,
\begin{align*}
\|\hat u-\tilde u\|_{L^p(0,T;X)}
\le c\left(\tau\sum_{n=1}^N
\|u^n-u^{n-1}\|_X^p\right)^{1/p}.
\end{align*}
For $n\ge2$, we have
\begin{align*}
u^n-u^{n-1}=\frac{2}{2-\alpha}(\tilde u^n-\tilde u^{n-1})-\frac{\alpha}{2-\alpha}(u^{n-1}-u^{n-2}).
\end{align*}
Since $\alpha/(2-\alpha)<1$, iteration of this identity and discrete Young's inequality yield
\begin{align*}
\Big(\tau\sum_{n=1}^N\|u^n-u^{n-1}\|_X^p\Big)^{1/p}
\le c&\Big(\tau\sum_{n=2}^N
\|\tilde u^n-\tilde u^{n-1}\|_X^p\Big)^{1/p}
+c\tau^{1/p}\|u^1-u^0\|_X.
\end{align*}
Moreover, using the identity $\tilde u^1-u^0=(1-\alpha/2)(u^1-u^0)$,
\begin{align*}
\Big(\tau\sum_{n=2}^N
\|\tilde u^n-\tilde u^{n-1}\|_X^p\Big)^{1/p}
&=\|\tilde u-\tilde u(\cdot-\tau)\|_{L^p(\tau,T;X)},\\
\tau^{1/p}\|u^1-u^0\|_X
&\le c\big(\|\tilde u\|_{L^p(0,\tau;X)}
+\tau^{1/p}\|u^0\|_X\big).
\end{align*}
Combining the preceding estimates proves the proposition.
\end{proof}

\section{Convergence for the backward subdiffusion}\label{cvg-bsp}
In this section, we extend the discrete maximal regularity framework to backward subdiffusion and establish the convergence of the corresponding PINN approximation. 

The backward subdiffusion problem reads: given the terminal data $ g \in Y \subset L^{2}(\Omega) $, find the initial data $ u_0 (x) $ such that
\begin{equation}\label{eq:backward-subdiffusion}
\left\{\begin{aligned}
    \partial_t^\alpha u - \Delta u = &\; f, \quad \mbox{in }\Omega \times (0,T],  \\
    u = &\; 0, \quad \mbox{on } \partial \Omega \times (0,T],  \\
    u (\cdot,T) = &\; g , \quad \mbox{in }\Omega . 
\end{aligned} \right.
\end{equation} 
The PINN loss for problem \eqref{eq:backward-subdiffusion} reads
\begin{align}
    \mathcal{G}_\ell^b (v) = &\; \| \bar{\partial}_{\tau}^\alpha v - \Delta \tilde{v} - \tilde{f} \|_{L^{p}(0,T;L^{2}(\Omega))}^{p} + \| \bar{\partial}_{\tau}^\alpha v - \Delta \tilde{v} -\tilde{f}\|_{\mathfrak{M}} + \| v(\cdot,T) - g \|_{Y}^{2} , \label{eq:loss-backward-subdiffusion}
\end{align} 
where the norm $ \| \cdot \|_{\mathfrak{M}} $ is defined by
\begin{equation}\label{M-space}
    \| \tilde{f} \|_{\mathfrak{M}} := \tau \sum_{n=1}^N t_{N-n+1}^{\frac{1}{p}-\frac{\alpha}{2} -1} \| \tilde{f}^n \|_{\dot H^1(\Omega)}.
\end{equation}
The additional term $\|\cdot\|_{\mathfrak{M}}$  plays the role of a regularization term and is crucial for controlling the ill-posedness of the backward problem.
The terminal data $ g $ belongs to the space $ Y := (\dot H^{2}(\Omega), \dot H^{4}(\Omega))_{1-\frac{1}{\alpha p},p} $.
Here, for $s \ge 0$, the space $\dot H^s (\Omega)$ is defined by
\begin{align*}
    \dot H^s (\Omega) := \Big\{ v \in L^2 (\Omega): \sum_{j=1}^\infty \lambda_j^s (v, \varphi_j)^2_\Omega < \infty \Big\},
\end{align*}
equipped with the norm
$ \| v \|_{\dot H^s (\Omega)}^2 := \sum_{j=1}^\infty \lambda_j^s (v, \varphi_j)^2_\Omega$,
where $0 < \lambda_1 < \lambda_2 \le \cdots \le \lambda_j \le  \cdots \to \infty$ are the eigenvalues of the Dirichlet Laplacian $-\Delta$, and $\{\varphi_j\}_{j=1}^\infty \subset H^2 (\Omega) \cap H_0^1 (\Omega)$ are the corresponding $L^2(\Omega)$ orthonormal eigenfunctions. 

Let $ \{u_\ell\}_\ell $ be a sequence of minimizers of the loss functional $ \mathcal{G}_\ell^b $ in the function class $ \mathcal{V}_\ell $. We impose the following assumption.
\begin{assumption}\label{assump-back}
The conditions in Assumption~\ref{assump0} hold with \eqref{eq:assump0} replaced by
    \begin{align}
        &\; \| \bar{\partial}_{\tau}^\alpha u \|_{L^{p}(0,T;L^{2}(\Omega))} + \| \Delta \tilde{u} \|_{L^{p}(0,T;L^{2}(\Omega))} \notag \\
        \le &\; c\big(\| \bar{\partial}_{\tau}^\alpha u - \Delta \tilde{u} \|_{L^{p}(0,T;L^{2}(\Omega))} +\| \bar{\partial}_{\tau}^\alpha u - \Delta \tilde{u} \|_{\mathfrak{M}}+ \| u_N \|_Y \big) , \label{eq:assump-back-0}
    \end{align} 
	for some $ 1 < \alpha p < 2 $.
\end{assumption}

By replacing Assumption~\ref{assump0} with Assumption~\ref{assump-back} in the proof of Theorem~\ref{main-thm-forward}, we obtain the following convergence result for backward subdiffusion. The proof is identical with that for Theorem~\ref{main-thm-forward}.
\begin{theorem}\label{main-thm-backward}
Let $ \{u_\ell\}_\ell $ be the minimizers of the loss \eqref{eq:loss-backward-subdiffusion} and Assumption~\ref{assump-back} holds. Then the piecewise linear temporal reconstructions $ \{ \hat{u}_\ell \}_\ell$ converge to the solution $ u $ of \eqref{eq:backward-subdiffusion} in $ L^{p}(0,T;X) $.
\end{theorem}

Below we verify that the BE, BDF2, L1, and fractional CN schemes satisfy \eqref{eq:assump-back-0}.

\subsection{BE scheme and BDF2 scheme}\label{subsec-BE-BDF2-backward} 
In this part, we prove that the BE and BDF2 schemes satisfy \eqref{eq:assump-back-0}. Let $ \{u^n\}_{n=0}^{N} $ satisfy the discrete equation
\begin{align}
    \bar{\partial}_{\tau}^\alpha u - \Delta \bar{u} = \bar{f},  \label{eq:BE-backward-subdiffusion-t}
\end{align} 
with the terminal condition $ u^N = g $. We split $ u^n = u_1^n + u_2^n $ such that $ \{u_1^n\} $ and $ \{u_2^n\} $ respectively satisfy
\begin{align}
    \bar{\partial}_{\tau}^\alpha u_1 - \Delta \bar{u}_1 &= 0,  \quad u_1^0 = u^0, \label{eq:BE-backward-subdiffusion-t1}\\
    \bar{\partial}_{\tau}^\alpha u_2 - \Delta \bar{u}_2 &= \bar{f}, \quad u_2^0 = 0. \label{eq:BE-backward-subdiffusion-t2}
\end{align}

We show the maximal regularity estimate for $ u_1 $ and $ u_2 $ separately. For $ u_1 $, the key step is to establish a stability estimate for $ u_1^0 $.
\begin{proposition}\label{est_u0XY}
    Let $ u_1 $ satisfy \eqref{eq:BE-backward-subdiffusion-t1}, with $ \bar{\partial}_{\tau}^\alpha $ being the BE or BDF2 scheme. Then
    \begin{align*}        
        \| u_1^0 \|_{X} \le c\| u_1^N \|_{Y} . %
    \end{align*}
This and the maximal regularity estimate \eqref{eq:assump0} for $ u_1 $ imply
    \begin{align*}
        \| \bar{\partial}_{\tau}^\alpha u_1 \|_{L^p(0,T;L^{2}(\Omega))} + \| \Delta \bar{u}_1 \|_{L^p(0,T;L^{2}(\Omega))} \le c\| u_1^N \|_{Y} .
    \end{align*} 
\end{proposition} 
\begin{proof}
By \cite[Section 4]{zhang2020numerical}, with the spectral decomposition, have
    \begin{align*}
        u_1^n = F_{\tau}^n u_1^0 = \sum_{j=1}^{\infty} F_{\tau}^n (\lambda_j) (u_1^0, \varphi_j) \varphi_j , \quad n = 1,2,\ldots,N,
    \end{align*} 
    where $ (\cdot,\cdot)$ denotes the $L^2(\Omega)$ inner product and $ F_{\tau}^n (\lambda) $ satisfies $F_\tau^0(\lambda)=1$ and
    \begin{align*}
        \bar{\partial}_{\tau}^\alpha F_{\tau}^n (\lambda) + \lambda F_{\tau}^n (\lambda) = 0, \quad n = 1,2,\ldots,N.
    \end{align*}
Note that $ F_{\tau}^n (\lambda) $ can be explicitly expressed as
    \begin{align*}
        F_{\tau}^n (\lambda) = \frac{1}{2 \pi i} \int_{\Gamma_{\theta,\sigma}^\tau} e^{z t_n} \delta_{\tau}^k (e^{-z \tau})^{\alpha - 1} \big( \delta_{\tau}^k (e^{-z \tau})^\alpha + \lambda \big)^{-1} \,\d z . %
    \end{align*}
By \cite[Corollary 4.1]{zhang2020numerical}, $ F_{\tau}^n (\lambda) $ satisfies 
    \begin{align}
      c_0(1+\lambda t_n^\alpha)^{-1} \leq   F_{\tau}^{n} (\lambda) \leq c_1({1+ \lambda t_n^\alpha})^{-1}. \label{eq:prop-11-2}
    \end{align}
For the $ L^{2}(\Omega) $ space, using the estimate \eqref{eq:prop-11-2}, we have
    \begin{align*}
        \| u_1^0 \|_{L^{2}(\Omega)}^2 = &\; \sum_{j=1}^{\infty} |(u_1^0, \varphi_j)|^2 \le c \sum_{j=1}^{\infty} (1 + \lambda_j T^\alpha)^2|F_{\tau}^N (\lambda_j) (u_1^0, \varphi_j)|^2  \notag \\
        \le &\; c \sum_{j=1}^{\infty} (1 + \lambda_j^{2}) |F_{\tau}^N (\lambda_j) (u_1^0, \varphi_j)|^2 = c \| u_1^N \|_{\dot{H}^{2}(\Omega)}^2 . 
    \end{align*}
For $ \dot H^{2}(\Omega) $ space, again using \eqref{eq:prop-11-2}, we have
    \begin{align}
        \| u_1^0 \|_{\dot H^{2}(\Omega)}^2 = &\; \sum_{j=1}^{\infty} \lambda_j^{2} |(u_1^0, \varphi_j)|^2 \le c\sum_{j=1}^{\infty} \lambda_j^{2}(1 + \lambda_j T^\alpha)^2|F_{\tau}^N (\lambda_j) (u_1^0, \varphi_j)|^2  \notag \\
        \le &\; c \sum_{j=1}^{\infty} (\lambda_j^{2} + \lambda_j^{4}) |F_{\tau}^N (\lambda_j) (u_1^0, \varphi_j)|^2 \le c \| u_1^N \|_{\dot{H}^{4}(\Omega)}^2 . \label{u_0-H0-esti}
    \end{align} 
	Hence, by real interpolation inequality as shown in \cite[Section 3.5]{bergh1976interpolation}, the preceding two estimates imply 
$\| u_1^0 \|_{X} \le c \| u_1^N \|_{Y}$.  
\end{proof} 

The maximal regularity estimate for $ u_2 $ follows directly from \eqref{eq:assump0} for the BE and BDF2 schemes. Combining the estimates for $ u_1 $ and $ u_2 $ gives
\begin{align}\label{eq:maximal-regularity-backward-1}
    &\; \| \bar{\partial}_{\tau}^\alpha u \|_{L^p(0,T;L^{2}(\Omega))} + \| \Delta \bar{u} \|_{L^p(0,T;L^{2}(\Omega))} \notag \\
    \le &\; \| \bar{\partial}_{\tau}^\alpha u_1 \|_{L^p(0,T;L^{2}(\Omega))} + \| \Delta \bar{u}_1 \|_{L^p(0,T;L^{2}(\Omega))} \notag \\
    &\; + \| \bar{\partial}_{\tau}^\alpha u_2 \|_{L^p(0,T;L^{2}(\Omega))} + \| \Delta \bar{u}_2 \|_{L^p(0,T;L^{2}(\Omega))} \notag \\
    \le &\; c \| u_1^N \|_{Y} + c \| \bar{f} \|_{L^p(0,T;L^{2}(\Omega))} \notag \\
    \le &\; c \| \bar{f} \|_{L^p(0,T;L^{2}(\Omega))} + c \| g \|_{Y} + c \| u_2^N \|_{Y}. 
\end{align} 

We have the following estimate on $ u_2^N $.
\begin{proposition}\label{lem-u2N-f}
For $ u_2 $ satisfying \eqref{eq:BE-backward-subdiffusion-t2}, we have 
$\| u_2^N \|_{Y} \le c\| \bar{f} \|_{\mathfrak{M}}$.
\end{proposition} 
\begin{proof}
Let
    \begin{align*}
        f (\xi) = \sum_{n=1}^\infty f^n \xi^n \quad\mbox{and}  \quad u_2 (\xi) = \sum_{n=1}^\infty u_2^n \xi^n .
    \end{align*} 
By the identity \eqref{eq:prop-7-2}, multiplying both sides of \eqref{eq:BE-backward-subdiffusion-t2} at step $ n $ by $\xi^n$ and summing over $n$ give
\begin{align*}
        u_2 (\xi) = &\; \big( \delta_{\tau}^k (\xi)^{\alpha} - \Delta \big)^{-1} f (\xi) %
    .%
\end{align*} 
Similar to the derivation of the identity \eqref{exp-barpartial-tau}, we have
\begin{align*}
        u_2^N = &\; \frac{\tau}{2 \pi i} \sum_{j=1}^N \int_{\Gamma_{\theta,\sigma}^\tau} e^{z t_{N-j}} \big( \delta_{\tau}^k (e^{-z \tau})^{\alpha} - \Delta \big)^{-1} f^j \, \d z %
        . %
\end{align*} 
By \cite[(1.14) and (1.19)]{jin2023numerical}, we have 
\begin{align*}
   \| \big( \delta_{\tau}^k (e^{-z \tau})^{\alpha} - \Delta \big)^{-1} f^j \|_{\dot H^{q}(\Omega)} \le C |z|^{-\frac{3-q}{2}\alpha} \| f^j \|_{\dot H^{1}(\Omega)}, \quad 1 \le q \le 3. %
\end{align*} 
By \cite[Theorem 3.5.3]{bergh1976interpolation}, $ (\dot H^{2}(\Omega), \dot H^{3}(\Omega))_{2-\frac{2}{\alpha p}, p} = (\dot H^{2}(\Omega), \dot H^{4}(\Omega))_{1-\frac{1}{\alpha p},p} $ for $ 1 < \alpha p < 2 $. By interpolation inequality, we obtain 
\begin{align*}
   \| \big( \delta_{\tau}^k (e^{-z \tau})^{\alpha} - \Delta \big)^{-1} f^j \|_{Y} \le C |z|^{\frac{\alpha}{2} - \frac{1}{p}} \| f^j \|_{\dot H^{1}(\Omega)}. %
\end{align*} 
Similarly to the argument for \eqref{eq:prop-7-10}, we obtain the desired estimate.
\end{proof} 

The estimate \eqref{eq:maximal-regularity-backward-1} and Proposition~\ref{lem-u2N-f} imply that BE and BDF2 schemes satisfy \eqref{eq:assump-back-0}.

\subsection{L1 scheme} 
In this part, we prove that the L1 scheme satisfies \eqref{eq:assump-back-0}. We split $ u^n = u_1^n + u_2^n $ such that $ \{u_1^n\} $ and $ \{u_2^n\} $ satisfy the schemes \eqref{eq:BE-backward-subdiffusion-t1} and \eqref{eq:BE-backward-subdiffusion-t2}, respectively. We prove two results analogous to Propositions~\ref{est_u0XY} and~\ref{lem-u2N-f}.  

First, to show the estimate for $ u_1^n$, we prove the stability of $ u_1^0 $ using the solution representation of $ u_1^n $. By \cite[Proposition 5.1]{jin2018discrete}, we have
\begin{align*}
    u_1^n = &\; u_1^0 + \frac{1}{2 \pi i} \int_{\Gamma_{\theta,\sigma}^\tau} e^{z t_{n+1}} \delta_{\tau}^1 (e^{-z \tau})^{-1} \Delta \big( \gamma_{\tau} (e^{-z \tau})^{\alpha} - \Delta \big)^{-1} u_1^0 \, \d z ,%
\end{align*} 
with $ \gamma_{\tau}(\xi) $ defined in \eqref{def_gamma_xi}.
Using the spectral decomposition $u_1^0 = \sum_{j=1}^{\infty} (u_1^0, \varphi_j)\varphi_j$ of $ u_1^0 $, 
 we have
\begin{align*}
    u_1^n = &\; \sum_{j=1}^{\infty} \mathcal{F}_{\tau}^n (\lambda_j)  (u_1^0, \varphi_j) \varphi_j , %
\end{align*}
with $ \mathcal{F}_{\tau}^n (\lambda_j) $ defined by
\begin{align*} \mathcal{F}_{\tau}^n (\lambda_j) = &\; 1 - \frac{\lambda_j}{2 \pi i} \int_{\Gamma_{\theta,\sigma}^\tau} e^{z t_{n+1}} \delta_{\tau}^1 (e^{-z \tau})^{-1} \big( \gamma_{\tau} (e^{-z \tau})^{\alpha} + \lambda_j \big)^{-1} \, \d z . %
\end{align*} 
Similar to the proof of \cite[Corollary 4.1]{zhang2020numerical}, Lemma~\ref{prop-generating-func-L1} yields 
\begin{align}
    c_0(1+\lambda_jt_n^\alpha)^{-1} \leq \mathcal{F}_{\tau}^n (\lambda_j) \leq c_1({1 + \lambda_j t_n^\alpha})^{-1} . \label{eq:prop-12-4}
\end{align} 
Using \eqref{eq:prop-12-4}, 
similar to the proof of Proposition~\ref{est_u0XY}, we obtain the following result.
\begin{proposition}\label{est_u0XY_L1}
    For $ u_1 $ satisfying \eqref{eq:BE-backward-subdiffusion-t1} with $ \bar{\partial}_{\tau}^\alpha $ being the L1 scheme, we have
\begin{align*}    \| u_1^0 \|_{X} \le c  \| u_1^N \|_{Y}.    \end{align*}
\end{proposition} 

Next, by \cite[Proposition 5.1]{jin2023numerical}, we have
\begin{align*}
    u_2^n = &\; \frac{\tau}{2 \pi i} \sum_{j=1}^n \int_{\Gamma_{\theta,\sigma}^\tau} e^{z t_{n-j}}  \big( \gamma_{\tau} (e^{-z \tau})^{\alpha} - \Delta \big)^{-1}  f^j \, \d z . %
\end{align*}

By Lemma~\ref{prop-generating-func-L1} and the proof of Proposition~\ref{lem-u2N-f}, we have the following bound for $ u_2^N $.
\begin{proposition}\label{lem-u2N-f-L1}
For $ u_2 $ satisfying \eqref{eq:BE-backward-subdiffusion-t2} with $ \bar{\partial}_{\tau}^\alpha$ being the L1 scheme, we have 
    \begin{align*}
        \| u_2^N \|_{Y} \le c \| \bar{f} \|_{\mathfrak{M}} .
    \end{align*} 
\end{proposition} 

\subsection{Fractional CN scheme} 
The derivation of \eqref{eq:assump-back-0} for the fractional CN scheme is similar to that for the BE and BDF2 schemes, and is  omitted.

\section{Numerical experiments and discussions}\label{numerics} 
Now we present numerical experiments to illustrate the theoretical findings. Consider the following subdiffusion problem:
\begin{align}\label{eq:subdiffusion-numerics}
	\partial_t^\alpha u (x,t) - \partial_{xx} u (x,t) = 0, \quad (x,t) \in (0,1) \times (0,T],
\end{align}
equipped with a zero Dirichlet boundary condition and initial condition $ u_0$. The exact solution $u$ to \eqref{eq:subdiffusion-numerics} is given by
\begin{align}\label{eq:subdiffusion-solution-numerics}
	u (x,t) = \sum_{n=1}^{\infty} 2 (u_0, \sin{(n \pi x)})_{L^2(0,1)} \sin{(n \pi x)} E_{\alpha}(-(n \pi)^{2}t^{\alpha}),
\end{align}
where $ E_{\alpha} (z) = \sum_{k=0}^{\infty} \frac{z^k}{\Gamma (\alpha k + 1)} $ is the Mittag-Leffler function \cite[Section 3.1]{jin2021fractional}. The reference solution can be computed by truncating the series expansion. 
We implement the PINN using a fully connected NN with $ 8 $ layers and $ 50 $ neurons in each hidden layer. The activation $\sigma$ is chosen to be hyperbolic tangent. We choose the BE scheme to discretize the fractional derivative $\partial_t^\alpha u$ in the loss. We enforce the zero Dirichlet boundary condition by multiplying the output of the NN by $ x(1-x) $. We use a deterministic sampling strategy to select $ N_t \times N_x $ uniform collocation points in the domain $Q_T$. The training of the PINN loss is performed using the Adam optimizer \cite{KingmaBa} with a learning rate $ 10^{-3} $. We measure the accuracy of the PINN using the $ L^p (0,T;H^{1}(0,1)) $ error between the predicted and reference solutions since the norm of $ X $ is inconvenient to compute. 

\subsection{Forward problem}
Consider the forward problem for $\alpha=0.25$, $0.5$, $0.75$, and choose $p=4$, $2$, $2$, respectively. The solutions are computed with $T=0.1$, $N_x=100$, and $N_t=10$. For $\alpha=0.25$ and $0.5$, we take $u_0 \in L^2(\Omega)$ such that
$u_0(x)=x(1-x)|x-0.5|^{-0.45}$,
which belongs to $L^2(\Omega)$ but not to $H^{0.05}(\Omega)$, whereas for $\alpha=0.75$, we take $u_0 = \sin(\pi x) \in H^1(\Omega)$, since the fractional Sobolev norm is inconvenient to evaluate numerically.

\begin{figure}[htbp]
	\centering
    \setlength{\tabcolsep}{0pt}
\begin{tabular}{ccc}	   \includegraphics[width=0.32\textwidth]{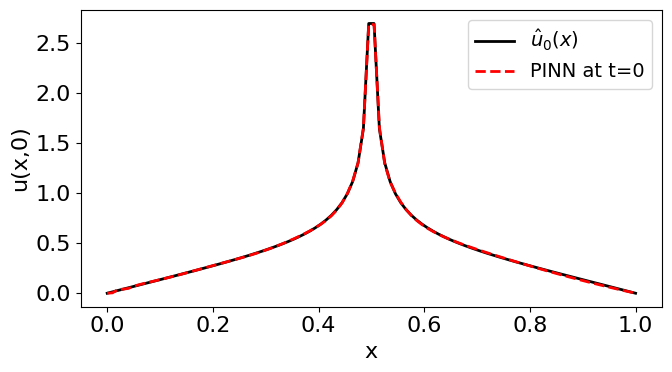}
&\includegraphics[width=0.32\textwidth]{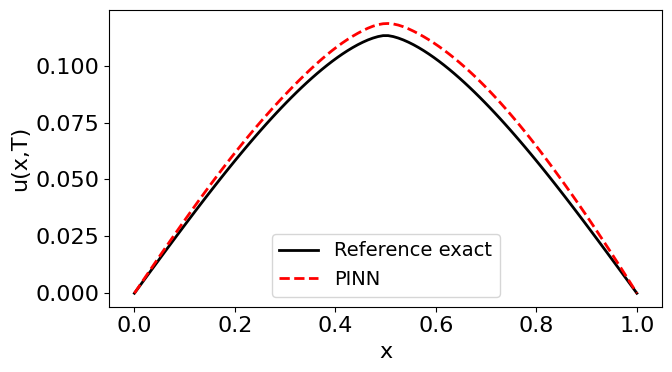}
&\includegraphics[width=0.32\textwidth]{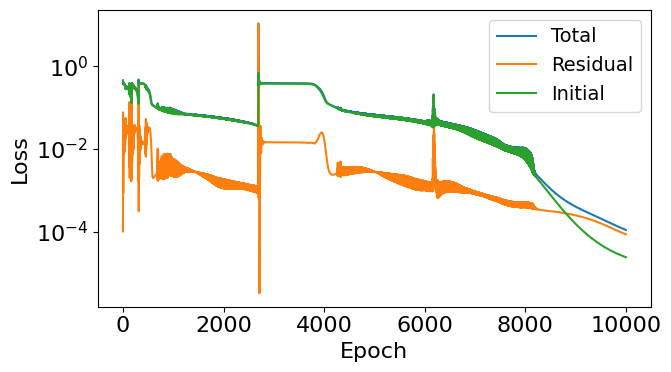}\\
\includegraphics[width=0.32\textwidth]{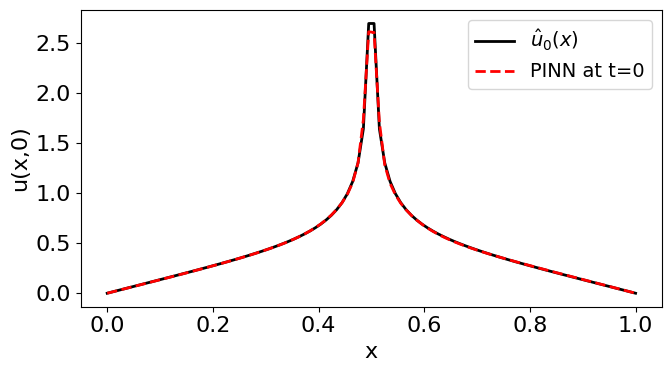}
&\includegraphics[width=0.32\textwidth]{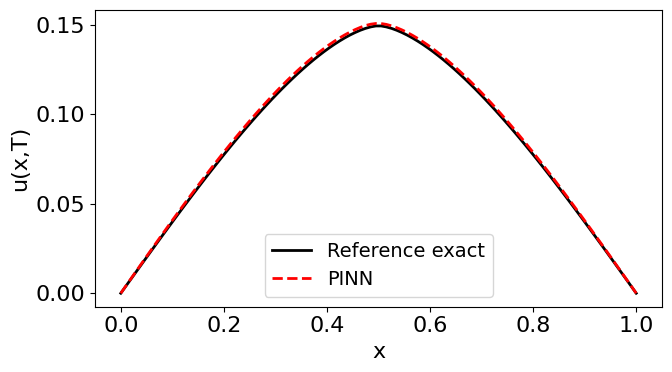}
&\includegraphics[width=0.32\textwidth]{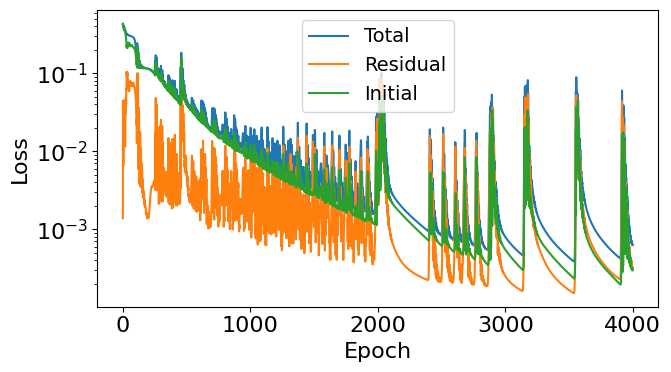}\\
\includegraphics[width=0.32\textwidth]{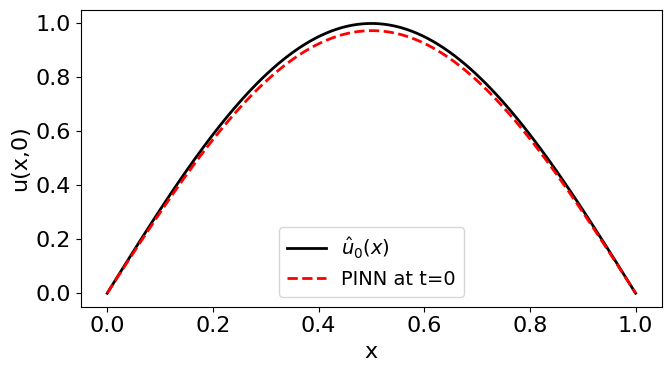}
&	\includegraphics[width=0.32\textwidth]{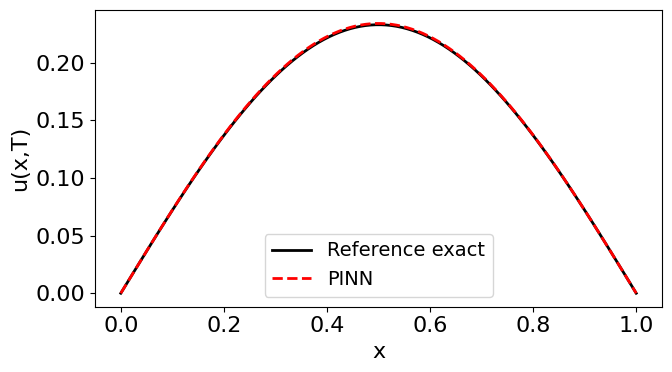}
&\includegraphics[width=0.32\textwidth]{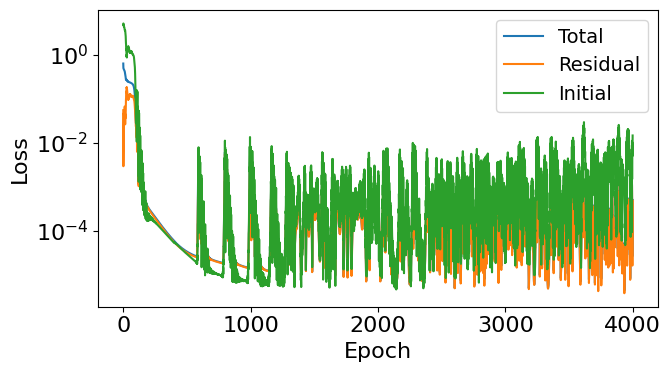}\\
$u_0$ & $u(\cdot,T)$ & training dynamics
\end{tabular}    
	\caption{The initial and terminal values, and training dynamics for the forward problem with $ \alpha = 0.25 $ (top), $\alpha=0.5$ (middle) and $\alpha=0.75$ (bottom).}
	\label{fig:forward25}
\end{figure}

The initial and terminal values, and the training dynamics for $ \alpha = 0.25$, 0.5, 0.75 are shown in Fig.~\ref{fig:forward25}. The PINN successfully captures the solution of the subdiffusion problem. The $ L^p (0,T;H^{1}(0,1)) $ errors are listed in Table~\ref{tab:forward-error}.

\begin{table}[htbp]
\caption{The $ L^p (0,T;H^{1}(0,1)) $ errors for the forward problem.}
\label{tab:forward-error}
\begin{center}
\begin{tabular}{c|ccc}
\hline 
$ \alpha $ & $ 0.25 $ & $ 0.5 $ & $ 0.75 $ \\
\hline
Values of $ p $ & $ 4 $ & $ 2 $ & $ 2 $ \\
$ L^p (0,T;H^{1}(0,1)) $ error & $ 3.058 \times 10^{-2} $ & $ 1.749 \times 10^{-2} $ & $ 1.050 \times 10^{-2} $ \\
\hline
\end{tabular}
\end{center} 

\end{table}

\subsection{Backward problem} 
\begin{figure}[htbp]
\centering\setlength{\tabcolsep}{0pt}
\begin{tabular}{ccc} \includegraphics[width=0.32\textwidth]{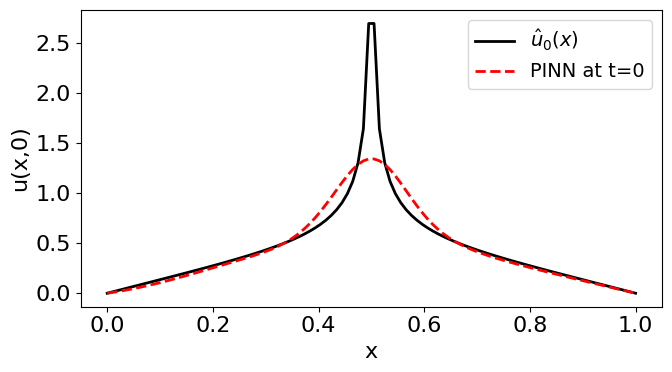}
&\includegraphics[width=0.32\textwidth]{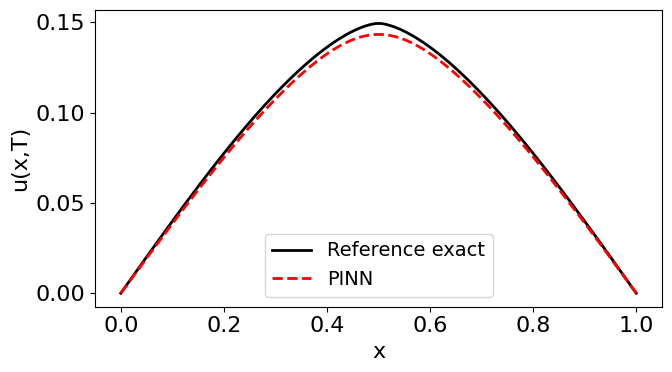}
&\includegraphics[width=0.32\textwidth]{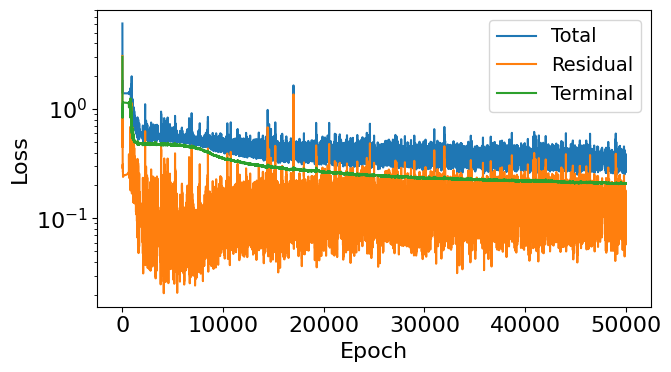}\\
\includegraphics[width=0.32\textwidth]{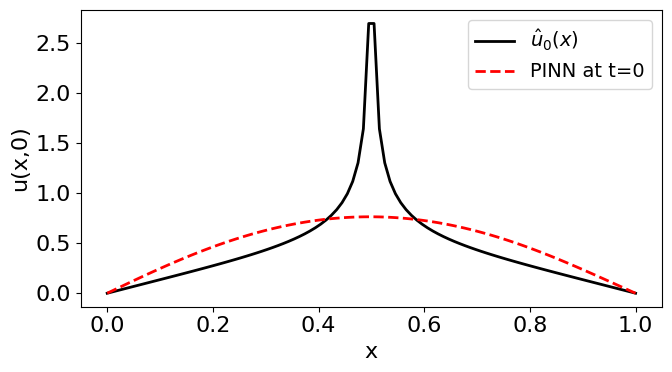}
&\includegraphics[width=0.32\textwidth]{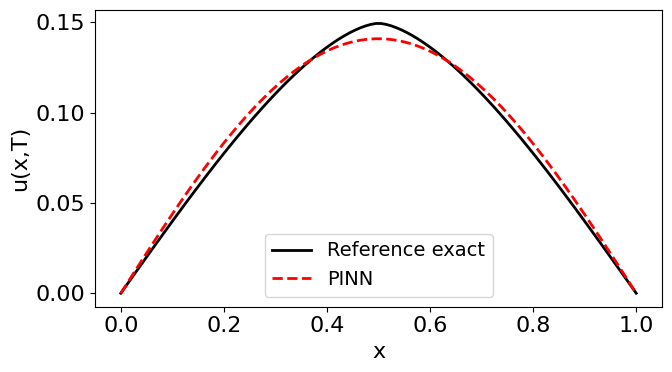}
&\includegraphics[width=0.32\textwidth]{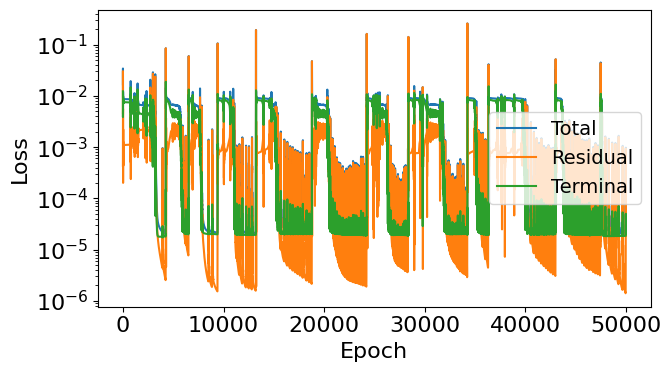}\\
$u_0$ & $u(\cdot,T)$ & training dynamics
\end{tabular}
\caption{The initial, terminal values and training loss for the backward problem using the regularized loss \eqref{eq:loss-backward-subdiffusion} (top) and the $L^2$ loss \eqref{eq:loss-bwd-subdiff-2} (bottom).}
	\label{fig:backward}
\end{figure}

Next, consider the backward subdiffusion problem. We choose the exact solution $u$ to be \eqref{eq:subdiffusion-solution-numerics} with the initial value $ u_0 (x) = x (1-x) | x-0.5 |^{-0.45} $. Then the terminal condition $g$ is taken to be $ u (\cdot,T) $. We take $ \alpha = 0.5 $, and $ p = 2 $. We compute the solutions with $ T = 0.1 $, $ N_x = 100 $, and $ N_t = 10 $. The initial value $u_0$, terminal value $u(\cdot,T)$ of the results and the training dynamics are shown in Fig.~\ref{fig:backward}. The $ L^2 (0,T;H^{1}(0,1)) $ error is about $ 1.377 \times 10^{-2} $. The results indicate that the PINN successfully reconstructs the solution of the backward problem.

We also test the PINN for the backward problem using only an $L^2$-type loss:
\begin{align}\label{eq:loss-bwd-subdiff-2}
	\mathcal{G}_\ell^b (u) = \| \bar{\partial}_{\tau}^\alpha u - \Delta \tilde{u} \|_{L^{2}(0,T;L^{2}(\Omega))}^2 + \| u_N - g \|_{L^{2}(\Omega)}^2 .
\end{align} 
The relevant numerical results are shown in Fig.~\ref{fig:backward}. The $ L^2 (0,T;H^{1}(0,1)) $ error is $ 2.569 \times 10^{-2} $. Then the PINN approximation fails to recover the initial data accurately. This observation confirms the necessity of incorporating suitable regularization terms into the loss functional for the backward problem, as indicated by the theoretical analysis and by the formulation \eqref{eq:loss-backward-subdiffusion}.

\section{Conclusion} \label{ccls}
In this paper, we have presented a rigorous convergence analysis of PINNs for linear subdiffusion and the associated backward problem. By formulating the PINN on the basis of classical time stepping schemes satisfying discrete maximal regularity, including the BE, BDF2, L1, and fractional CN schemes, we have established a theoretical framework that connects deep learning-based solvers with classical numerical analysis for time-fractional evolution equations.

Our main theoretical result shows that the minimizers of the PINN loss converge to the exact solution that belongs to the critical space. The analysis is based on $\liminf$ and $\limsup$ arguments, Aubin--Lions compactness lemma, and discrete integration-by-parts formulas for fractional derivatives. We have further verified that several commonly used time-stepping schemes satisfy the requisite discrete maximal regularity and consistency properties. In addition, the theory was extended to backward subdiffusion, for which convergence can still be established despite its ill-posedness, provided that suitable regularization is incorporated into the loss.
Both the theoretical analysis and the numerical experiments indicate that an $L^2$-type loss alone is generally insufficient to ensure convergence of the PINN approximation, and that additional regularization terms are essential, especially for the backward problem. This observation highlights the importance of combining approximation, stability, and regularization in the design of PINNs  for fractional PDEs. Overall, this work provides a first rigorous theoretical foundation for the application of PINNs to time fractional diffusion and contributes to a deeper understanding of learning-based methods.

Several interesting directions await future research. One important extension is to establish analogous convergence results for nonlinear subdiffusion and more general coupled systems. It is also of interest to investigate the effect of noisy sparse observations, particularly in the practical backward problem. It would also be of interest to study adaptive time-stepping strategies and optimization procedures tailored to the nonlocal structure of fractional derivatives, with the goal of improving both computational efficiency and accuracy.

\bibliographystyle{abbrv}
\bibliography{ref}

\end{document}